\documentclass[aap,seceqn,citesort,MSNbibl,dvips]{arximspdf}

\doi{10.1214/09-AAP610}
\volume{19}
\issue{6}
\pubyear{2009}
\firstpage{2270}
\lastpage{2300}

\makeatletter

\newcommand{\iint}{\int\!\!\int}

\newproclaim{condition}{Condition}[section]
\newtheorem{theorem}[condition]{Theorem}
\newtheorem{lemma}[condition]{Lemma}
\newtheorem{prop}[condition]{Proposition}

\newproclaim{Remark}{Remark}

\makeatother

\begin{document}
\begin{frontmatter}

\title{Limit theorems for additive functionals of~a~Markov chain}
\runtitle{Stable limit laws for Markov chains}

\begin{aug}
\author[A]{\fnms{Milton} \snm{Jara}\thanksref{t1}\ead[label=e1]{jara@ceremade.dauphine.fr}},
\author[B]{\fnms{Tomasz} \snm{Komorowski}\thanksref{t2,t3}\ead[label=e2]{komorow@hektor.umcs.lublin.pl}\ead
[label=u1,url]{http://hektor.umcs.lublin.pl/\texttildelow komorow}} and
\author[A]{\fnms{Stefano} \snm{Olla}\corref{}\thanksref{t3,t4}\ead[label=e3]{olla@ceremade.dauphine.fr}\ead
[label=u3,url]{http://www.ceremade.dauphine.fr/\texttildelow olla}}
\runauthor{M. Jara, T. Komorowski and S. Olla}
\affiliation{CEREMADE Universit\'{e} de Paris Dauphine,
Institute~of~Mathematics, UMCS~and~CEREMADE Universit\'{e} de Paris Dauphine}
\address[A]{M. Jara\\
S. Olla\\
CEREMADE, Universit\'{e} de Paris Dauphine\\
Place du Mar\'{e}chal De Lattre De Tassigny\\
75775 Paris Cedex 16\\
France\\
\printead{e1}\\
\phantom{E-mail: }\printead*{e3}\\
\printead{u3}}
\address[B]{T. Komorowski\\
Institute of Mathematics, UMCS\\
pl. Marii Curie-Sk\l odowskiej 1\\
Lublin 20-031\\
Poland\\
\printead{e2}\\
\printead{u1}}
\end{aug}

\thankstext{t1}{Supported by the Belgian Interuniversity Attraction
Poles Program P6/02, through the network NOSY (Nonlinear systems,
stochastic processes and statistical mechanics).}
\thankstext{t2}{Supported in part by Polish MNiSW Grant 20104531.}
\thankstext{t3}{Supported in part by EC FP6 Marie Curie ToK programme SPADE2,
MTKD-CT-2004-014508 and Polish MNiSW SPB-M.}
\thankstext{t4}{Supported by French ANR LHMSHE BLAN07-2184264.}

\received{\smonth{10} \syear{2008}}
\revised{\smonth{4} \syear{2009}}

%
\begin{abstract}
Consider a Markov chain $\{X_n\}_{n\ge0}$ with an ergodic probability
measure~$\pi$. Let $\Psi$ be a function on the state space of the chain,
with $\alpha$-tails with respect to $\pi$, $\alpha\in(0,2)$.
We find sufficient
conditions on the probability transition to prove convergence in law
of $N^{1/\alpha}\sum_n^N \Psi(X_n)$ to an $\alpha$-stable
law. A ``martingale approximation'' approach and a ``coupling''
approach
give two different sets of conditions. We extend these results to
continuous time Markov jump processes $X_t$, whose skeleton chain
satisfies
our assumptions. If waiting times between jumps have finite expectation,
we prove convergence of $N^{-1/\alpha}\int_0^{Nt} V(X_s) \,ds$ to a
stable process.
The result is applied to show that an appropriately scaled limit of
solutions of a
linear Boltzman equation is a solution of the fractional diffusion equation.
\end{abstract}

%
\begin{keyword}[class=AMS]
\kwd[Primary ]{60F05}
\kwd{60F17}
\kwd[; secondary ]{76P05}.
\end{keyword}
\begin{keyword}
\kwd{Stable laws}
\kwd{self-similar L\'{e}vy processes}
\kwd{limit theorems}
\kwd{linear Boltzmann equation}
\kwd{fractional heat equation}
\kwd{anomalous heat transport}.
\end{keyword}

\pdfkeywords{60F05, 60F17, 76P05, Stable laws,
self-similar Levy processes,
limit theorems,
linear Boltzmann equation,
fractional heat equation,
anomalous heat transport,}

\end{frontmatter}

\section{Introduction}

Superdiffusive transport of energy is generically observed in a certain
class of one-dimensional systems. This can be seen numerically in
chains of anharmonic
oscillators of the Fermi--Pasta--Ulam type and experimentally in carbon
nanotubes (see \cite{sll} for a physical review). The nature of the stochastic
processes describing these
emerging macroscopic behaviors is a subject of a vivid debate in the
physical literature and remarkably few mathematical results are present for
deterministic microscopic models.

The macroscopic behavior of the energy in
a chain of harmonic oscillators with the Hamiltonian dynamics
perturbed by stochastic terms conserving energy and momentum has been
studied in \cite{basollaspohn}. The density of energy distribution
over spatial and momentum variables, obtained there in a proper
kinetic limit, satisfies \textit{a linear phonon Boltzmann equation},
%
%
\begin{equation}
\label{eq:boltz}\hspace*{29pt}
\partial_t u(t,x,k)+\omega'(k)\,\partial_x u(t,x,k)= \int R(k,k')
\bigl(u(t,x,k')- u(t,x,k) \bigr) \,dk'.
\end{equation}
As we have already mentioned, $u(t,x,k)$ is the density at time $t$ of
energy of waves of
Fourier's mode $k\in[0,1]$, and the velocity $\omega'(k)$ is the
derivative of the dispersion relation of the lattice.

We remark at this point that (\ref{eq:boltz}) appears also as
a limit of scaled wave, or Schr\"{o}dinger equations in a random medium
with fast oscillating coefficients and initial data. It is sometimes
called, in that context, \textit{the radiative transport equation} (see,
e.g., \cite{BPRNonlin,ErdosYau,fannjiang,LukSpohn,spohn}, or
monography \cite{foque} for more details on this subject).

Since the kernel $R(k,k')$ appearing in (\ref{eq:boltz}) is positive,
this equation has an easy probabilistic interpretation as a
forward equation for the evolution of the density of a Markov
process $(Y(t), K(t))$ on $\mathbb R\times[0,1]$. In fact, here
$K(t)$ is an autonomous jump process on $[0,1]$ with jump rate
$R(k,k')$, and $Y(t)= \int_0^t \omega'(K(s)) \,ds$ is an
additive functional of $K(t)$. Momentum conservation in the
microscopic model imposes a very slow jump rate for small $k$:
$R(k,k')\sim k^2$ as $k\sim0$, while velocity $\omega'(k)$
remains of order 1 even for small $k$.
So when $K(t)$ has a small value, it may stay unchanged for a long
time, as does the velocity of $Y(t)$. This is the mechanism that
generates on a macroscopic scale the superdiffusive behavior of $Y(t)$.

The above example has motivated us to study
the following general question.
Consider a Markov chain $\{X_n, n\ge0\}$ taking values in a
general Polish metric space $(E,$d$)$. Suppose that $\pi$ is a
stationary and ergodic probability Borel measure for this chain. Consider
a function $\Psi\dvtx E\to\mathbb R$ and
$S_N:=\sum_{n=0}^{N-1}\Psi(X_n)$. If $\Psi$ is centered with respect
to $\pi$, and possesses a second moment, one expects that the central
limit theorem holds for $N^{-1/2}S_N$, as $N\to+\infty$. This, of
course, requires some assumptions about the rate of the decay of
correlations of the chain, as well as hypotheses about its dynamics.
If $\Psi$ has an infinite second moment and its tails satisfy a power
law, then one expects, again under
some assumption on the transition probabilities, convergence of the
laws of $N^{-1/\alpha}S_N$, for an appropriate $\alpha$ to the
corresponding stable law.

In 1937 W. Doeblin himself looked at this natural question in
his seminal article~\cite{doeblin}. In the final lines of this paper,
he observes that the method of dividing the sum into independent
blocks, used in the paper to show the central limit theorem for
countable Markov chains, can be used also in the infinite variance
situation. A~more complete proof, along the line of Doeblin's idea,
can be found in an early paper of Nagaev
\cite{nagaev}, assuming a strong \textit{Doeblin} condition.

Starting from the early sixties, another, more analytical approach,
has been developed for proving central limit theorems for Markov
chains, based on a martingale approximation of the additive functional.
By solving (or by approximating the solution of) the Poisson equation
$(I-P) u = \Psi$ 
where $P$ is the transition probability matrix,
one can decompose the sum $S_N$ into a martingale plus a negligible
term, thus
reducing the problem to a central limit theorem for martingales.
This is exploited by Gordin (see \cite{gordin}) when $P$ has a
spectral gap. In the
following decades, much progress has been achieved using this approach.
It has found applications in stochastic homogenization, random walks in
random environments and interacting
particle systems (i.e., infinite-dimensional problems, where
renewal arguments cannot be applied),
culminating in the seminal paper of
Kipnis and Varadhan \cite{kipnisvaradhan} where reversibility of the
chain is
exploited in an optimal way (see also
\cite{demassi,goldstein,derrenniclin}). For nonreversible chains
there are still open problems (see \cite{maxwellwoodrofe}
and the review paper \cite{olla} for a more detailed list).

As far as we know, the martingale approximation approach has not been
developed in the case of convergence to stable laws of functionals of
Markov chains, even though
corresponding theorems of martingale convergence have been
available for a while (cf., e.g., \cite{durrettresnick,browneagleson}).
The present article is a first step in this direction.

More precisely, we are concerned with the limiting behavior of
functionals formed over functions $\Psi$ with heavy tails that satisfy a
power law decay, that is, $\pi(\Psi>\lambda)\sim c_*^+\lambda
^{-\alpha}$
and $\pi(\Psi<-\lambda)\sim c_*^-\lambda^{-\alpha}$ for
$\lambda\gg1$ with $\alpha\in(0,2)$.
We prove sufficient conditions under which
the laws of the functionals of the form $N^{-1/\alpha}S_N$ converge
weakly to $\alpha$-stable laws, as $N\to+\infty$.
Theorem
\ref{stable-discrete-1}
is
proven by
martingale approximation, under a spectral gap condition.

We also give a proof by a more \textit{classical} renewal
method based on a coupling technique inspired by
\cite{cc}.
The coupling argument gives a simpler proof
but under more restrictive assumptions on the form of the
probability transition (cf. Condition \ref{doeblin}). We point
out, however, that such hypotheses are of local nature, in the sense
that they involve only the behavior of the process around the
singularity. In particular, the spectral gap condition (which is a
global condition) can be relaxed in this coupling approach,
to a moment bound for some regeneration times associated to the
process (cf. Theorem \ref{thm-coupl}).

Next, we apply these results to a continuous time Markov jump process
$\{X_t, t\ge0\}$ whose skeleton chain satisfies the assumptions made
in the respective parts of Theorem \ref{stable-discrete-1}. We prove
that if the mean waiting
time $t(x)$ has a finite moment with respect to the invariant measure
$\pi$ and the tails of $V(x)t(x)$ obey the power laws, as above, then
finite-dimensional distributions of
the scaled functional of the form $N^{-1/\alpha}\int_0^{Nt}V(X_s)\,ds$
converge to the
respective finite-dimensional distribution of a stable process (see
Theorem \ref{thm-main-3}).

Finally, these results are applied to deal with the limiting
behavior of the solution $u(t,x,k)$
of the linear Boltzmann equation (\ref{eq:boltz})
in the spatial dimension $d=1$.
We prove that the long-time, large-scale limit
of solutions of such an equation converges to the solution of the
fractional heat equation
\[
\partial_t \bar u(t,x) = -(-\partial_x^2)^{3/4} \bar u(t,x) ,
\]
corresponding to a stable process with exponent $\alpha=3/2$. Both
approaches (i.e., martingale approximation and coupling) apply to this
example.

\textit{Note added to the second version}: After completing the first version
of the present paper \cite{v1},
we have received a preprint by Mellet, Mischler and Mouhot \cite{mmm} that
contains a completely analytical proof of the convergence of the
solution of a linear Boltzmann equation to a fractional diffusion. The
conditions assumed in \cite{mmm} imply the same spectral gap
condition as in our Theorem \ref{stable-discrete-1}; consequently the
corresponding result in \cite{mmm} is related to our Theorem \ref{thm-main-3}.

\section{Preliminaries and statements of the main results}

\subsection{Some preliminaries on stable laws}

In this paper we shall consider three types of stable laws.
When $\alpha\in(0,1)$, we say that $X$ is distributed according to a
stable law
of type I if its characteristic function is of the form
${\mathbb E}e^{i\xi X}=e^{\psi(\xi)}$, where the L\'{e}vy exponent equals
%
%
\begin{equation}
\label{type-1}
\psi(\xi):=\alpha\int_{{\mathbb R}}(e^{i\lambda
\xi}-1)|\lambda|^{-1-\alpha}c_*(\lambda)\,d\lambda
\end{equation}
and
%
%
\begin{equation}
\label{c-levy}
c_*(\lambda):=\cases{
c_*^-,&\quad when $\lambda<0$,\vspace*{2pt}\cr
c_*^+,&\quad when $\lambda>0$,}
\end{equation}
where $c_*^-,c_*^+\ge0$ and $c_*^-+c_*^+>0$.
The stable law is of type II
if $\alpha\in(1,2)$ and its L\'{e}vy exponent equals
%
%
\begin{equation}
\label{type-2}
\psi(\xi):=\alpha\int_{{\mathbb R}}(e^{i\lambda\xi}-1-i\lambda
\xi
)|\lambda|^{-1-\alpha}c_*(\lambda)\,d\lambda.
\end{equation}
Finally, the stable law is of type III is $\alpha=1$ and its
L\'{e}vy exponent
equals
%
%
\begin{eqnarray}
\label{eq:75d}
&&\psi(\xi):=
\int_{{\mathbb R}}\bigl(e^{i\xi\lambda}-1-i\xi\lambda1_{[-1,1]}(\lambda
)\bigr)|\lambda
|^{-2}c_*(\lambda)\,d\lambda.
\end{eqnarray}

We say that $\{Z(t), t\ge0\}$ is a stable process of type I (resp.,
II, or III) if \mbox{$Z(0)=0$} and it is a
process with independent increments such that $Z(1)$ is distributed
according to
a stable law of type I (resp., II, or III).


\subsection{A Markov chain}
\label{markov-chain}

Let $(E,$d$)$ be a Polish metric space, $\mathcal{E}$ its Borel
$\sigma$-algebra. 
Assume that $\{X_n, n\ge0\}$ is a Markov chain with the state space
$E$ and $\pi$---the law of $X_0$---is an \textit{invariant} and \textit{ergodic} measure for the chain. Denote by $P$ the transition operator
corresponding to the chain. Since $\pi$ is invariant it can be
defined, as a positivity preserving linear contraction, on any $L^p(\pi)$ space for $p\in[1,+\infty]$.
\begin{condition}\label{power}
Suppose that $\Psi\dvtx E\to{\mathbb R}$ is Borel measurable such that
there exist $\alpha\in(0,2)$ and two
constants $ c_*^+, c_*^-$ satisfying $c_*^++c_*^->0$ and
%
%
\begin{eqnarray}
\label{eq:tails}
\lim_{\lambda\to+\infty} \lambda^\alpha\pi(\Psi\ge\lambda) &=& c_*^+,
\nonumber\\[-8pt]\\[-8pt]
\lim_{\lambda\to+\infty} \lambda^\alpha\pi(\Psi\le-\lambda)
&=& c_*^-.\nonumber
\end{eqnarray}
\end{condition}

Condition (\ref{eq:tails}) guarantees that
$\Psi\in L^\beta(\pi)$ for any $\beta< \alpha$.

In the case of $\alpha\in(1,2)$, we will always assume that
$\int\Psi \,d\pi= 0$.
We are interested in the asymptotic behavior of
$S_N:=\sum_{n=1}^N \Psi(X_n)$. We are looking for sufficient
conditions on the chain, which guarantee that the laws of $N^{-1/\alpha
} S_N$ converge
to a $\alpha$-stable law, as $N\to+\infty$.

We present two different approaches (by martingale approximation and by
coupling) with two separate set of conditions.

\subsection{The martingale approach result}
\label{sec:martingale-approach}

We suppose that the chain satisfies:
%

\begin{condition}[(Spectral Gap Condition)]\label{sg}
%
%
\begin{equation}
\label{eq:8}
\sup\bigl[\|Pf\|_{L^2(\pi)}\dvtx f\perp1, \|f\|_{L^2(\pi)}=1\bigr]=a<1.
\end{equation}
\end{condition}

Since $P$ is also a contraction in $L^1(\pi)$ and $L^\infty(\pi)$ we
conclude, via the Riesz--Thorin interpolation theorem, that for any
$p\in[1,+\infty)$,
%
%
\begin{equation}
\label{eq:9}
\|Pf\|_{L^p(\pi)}\le a^{1-|2/p-1|}\|f\|_{L^p(\pi)},
\end{equation}
for all $f\in
L^p(\pi)$, such that $\int f\,d\pi=0$.


In addition, we assume that the tails of $\Psi$ under the invariant
measure do not differ very much from those with respect to the
transition probabilities. Namely, we suppose:
\begin{condition}
\label{decomp}
There exists a measurable family of Borel measures $Q(x, dy)$ and a
measurable, nonnegative function $p(x,y)$ such that
%
%
\begin{eqnarray}\label{transition}
P(x,dy)&=&p(x,y)\pi(dy)+Q(x,dy)\qquad \mbox{for all }x\in E,
\\
\label{eq:4}
C(2):\!&=& \sup_{y\in E} \int p^2(x,y)\pi(dx)<+\infty
\end{eqnarray}
and
%
%
\begin{equation}
\label{new}
Q(x, |\Psi|\ge\lambda) \le
C\int_{[|\Psi(y)|\ge\lambda]} p(x,y)\pi(dy)\qquad\forall x\in
E,\lambda\ge0.
\end{equation}
\end{condition}

A simple consequence of (\ref{transition}) and the fact that $\pi$ is
invariant is that
%
%
\begin{equation}\label{eq:2}
\int p(x,y)\pi(dy)\le1 \quad\mbox{and}\quad\int p(y,x)\pi(dy)\le
1\qquad \forall x\in E.
\end{equation}

If $\alpha\in(1,2)$
then, in particular, $\Psi$ possesses the first absolute
moment.
%
\begin{theorem}\label{stable-discrete-1}
We assume here Conditions \ref{power}--\ref{decomp}.
\begin{longlist}
\item Suppose $\alpha\in(1,2)$, $\Psi$ is centered.
Furthermore, assume that
for some \mbox{$\alpha'>\alpha$}, we have
%
%
\begin{equation}
\label{020411}
\|P\Psi\|_{L^{\alpha'}(\pi)}<+\infty.
\end{equation}
Then the law of $N^{-1/\alpha}S_N$ converges weakly,
as $N\to+\infty$, to a stable law of type~\textup{II}.
\item If $\alpha\in(0,1)$, then
the law of $N^{-1/\alpha}S_N$ converges weakly,
as $N\to+\infty$, to a stable law of type \textup{I}.
\item When $\alpha=1$, assume that for some $\alpha'>1$, we have
%
%
\begin{equation}
\label{020411b}
{\sup_{N\ge1}}\|P\Psi_N\|_{L^{\alpha'}(\pi)}<+\infty,
\end{equation}
where $\Psi_N:=\Psi1[|\Psi|\le N]$. Let $c_N:=\int\Psi_N\,d\pi$.
Then,
the law of $N^{-1}(S_N-Nc_N)$ converges weakly,
as $N\to+\infty$, to a stable law of type \textup{III}.
\end{longlist}
\end{theorem}
\begin{Remark*}
A simple calculation shows that in case (iii) $c_N=(c+o(1))\log N$ for
some constant $c$.
\end{Remark*}

\subsection{The coupling approach results}
\label{sec:coupl-appr-results}
\begin{condition}
\label{doeblin}
There exists a measurable function $\theta\dvtx E \to[0,1]$,
a probability $q$ and a
transition probability $Q_1(x, dy)$, such that
\[
P(x, dy) = \theta(x) q(dy) + \bigl(1-\theta(x)\bigr) Q_1(x, dy).
\]
Furthermore, we assume that
%
%
\begin{equation}
\label{eq:24}
\bar\theta:= \int\theta(x) \pi(dx) > 0
\end{equation}
and that
the tails of distribution of $\Psi$ with respect to $Q_1(x, dy)$ are
uniformly lighter than its tails with respect to $q$,
%
%
\begin{equation}\label{tailsQ1}
\lim_{\lambda\to\infty} \sup_{x\in E} \frac{Q_1(x, |\Psi| \ge
\lambda)}{q( |\Psi| \ge\lambda)} = 0.
\end{equation}
\end{condition}

Clearly, because of (\ref{tailsQ1}),
the function $\Psi$ satisfies condition (\ref{eq:tails})
also with respect to the measure $q$, but with different constants.
%
%
\begin{eqnarray}
\label{eq:tails1}
\lim_{\lambda\to+\infty} \lambda^\alpha q(\Psi>\lambda) &=& c_*^+
\bar\theta^{-1},
\nonumber\\[-8pt]\\[-8pt]
\lim_{\lambda\to+\infty} \lambda^\alpha q(\Psi< -\lambda) &=&
c_*^- \bar\theta^{-1}.\nonumber
\end{eqnarray}

The purpose of Condition \ref{doeblin} is that it permits to define
a Markov chain $\{(X_n, \delta_n), n\ge0\}$ on $E\times\{0,1\}$
such that
%
%
\begin{eqnarray}
\label{eq:7}
\mathbb P(\delta_{n+1} =0| X_n =x,\delta_n=\epsilon) &=& \theta(x),
\nonumber\\
\mathbb P(\delta_{n+1} =1| X_n =x,\delta_n=\epsilon) &=& 1-\theta
(x), \nonumber\\[-8pt]\\[-8pt]
\mathbb P(X_{n+1}\in A| \delta_{n+1} =0, X_n = x,\delta_n=\epsilon)
&=& q(A) ,\nonumber\\
\mathbb P(X_{n+1}\in A| \delta_{n+1} =1, X_n = x,\delta_n=\epsilon) &=&
Q_1(x,A)\nonumber
\end{eqnarray}
for $n\ge0$. We call this Markov chain the \textit{basic coupling}. It
is clear that
the marginal chain $\{X_n, n\ge0\}$ has probability transition $P$.
The dynamics of $\{(X_n,\delta_n), n\ge0\}$ are easy to understand. When
$X_n=x$, we choose $X_{n+1}$ according to the distribution $q(dy)$
with probability $\theta(x)$, and according to the distribution
$Q_1(x,dy)$ with probability $1-\theta(x)$.

Let $\kappa_n$ be the $n$th zero in the sequence $\{\delta_n, n\ge
0\}$. In a
more precise way, define $\kappa_0:=0$, and for $i \geq1$,
\[
\kappa_i: = \inf\{n> \kappa_{i-1}, \delta_n =0\}.
\]
%
Notice that the sequence $\{\kappa_{i+1}-\kappa_i, i
\geq1\}$ is i.i.d., and $\mathbb E(\kappa_{i+1}-\kappa_i)= \bar
\theta^{-1}$. We call the sequence $\{\kappa_n, n\ge1\}$ the \textit{regeneration times}. 

Observe that, for any $i \geq1$, the distribution of $X_{\kappa_i}$ is
given by $q(dy)$. In particular, $X_{\kappa_i}$ \textit{is independent
of} $\{X_0,\ldots,X_{\kappa_i-1}\}$. Therefore, the blocks
\[
\{(X_{\kappa_{i}},\delta_{\kappa_{i}}),\ldots,(X_{\kappa_{i+1}-1},
\delta_{\kappa_{i+1}-1})\}
\]
are independent. The dynamics for each one of these blocks is easy to
understand. Start a Markov chain $\{X_n^1, n\ge0\}$ with initial
distribution $q(dy)$ and transition probability $Q_1(x,dy)$. At each
step $n$, we stop the chain with probability $\theta(X_n^1)$. We call
the corresponding stopping time $\kappa_1$. Each one
of the blocks, except for the first one, has a distribution
$\{(X_{0}^1,0),(X_{1}^1,1),\ldots,\break (X_{\kappa_1-1}^1,1)\}$. The first
block is constructed in the same way, but
starts from $X_0^1=X_0$ instead of with the law $q(dy)$.
Now we are ready to state:
\begin{condition}
\label{regen}
\[
\sum_{n =1}^\infty n^{1+\alpha}\sup_x {\mathbb P}(\kappa_1\geq n|X_0=x)
<+\infty.
\]
\end{condition}
\begin{theorem}\label{thm-coupl}
Suppose that $\alpha\in(1,2)$ and $\Psi$ is centered under $\pi$,
or $\alpha\in(0,1)$.
Then under Conditions \ref{power}, \ref{doeblin} and \ref{regen},
the law of
$N^{-1/\alpha}S_N$ converges to an $\alpha$-stable law.
\end{theorem}
%

\subsection{An additive functional of a continuous time jump process}

Suppose that $\{\tau_n, n\ge0\}$ is a sequence of i.i.d. random
variables, independent of
$\mathcal{F}:=\sigma(X_0,X_1,\ldots)$ and such that $\tau_0$ has
exponential distribution with parameter $1$. Suppose that $t\dvtx E\to
(0,+\infty)$ is a measurable function such that $t(x)\ge t_*>0$, $x\in
E$. Let
%
%
\begin{equation}
\label{t-N}
t_N:=\sum_{n=0}^Nt(X_n)\tau_n.
\end{equation}
One can define a compound Poisson process $X_t=X_n$, $t\in[t_N,t_{N+1})$.
It is Markovian; see, for example, Section 2 of Appendix 1, pages
314--321, of \cite{kipnislandim} with the generator
%
%
\begin{equation}
\label{gen-jump-proc}
Lf(x)=t^{-1}(x)\int[f(y)-f(x)] P(x,dy),\qquad f\in B_b(E).
\end{equation}
Here $B_b(E)$ is the space of bounded and Borel measurable functions on $E$.
Let
%
%
\begin{equation}\label{bar-t}
\bar t:=\int t\,d\pi<+\infty.
\end{equation}
%
Suppose $V\dvtx E\to\mathbb R$ is measurable and
$\Psi(x):=V(x)t(x)$ satisfies condition (\ref{eq:tails}).
We shall be concerned with the limit of the scaled processes,
%
%
\begin{equation}
\label{eq:101}
Y_N(t):=\frac{1}{N^{1/\alpha}}\int_0^{Nt}V(X(s))\,ds,\qquad t\ge0,
\end{equation}
as $N\to+\infty$. 
Then $\bar t^{-1}t(x)\pi(dx)$ is an ergodic, invariant probability
measure for $\{X_t, t\ge0\}$.
Our result can be formulated as follows.
\begin{theorem}
\label{thm-main-3}
\textup{(i)} Suppose that $\alpha\in(1,2)$ and that the assumptions of either
part \textup{(i)} Theorem \ref{stable-discrete-1}, or of Theorem \ref{thm-coupl},
hold. Then, the convergence of finite-dimensional distributions takes
place to a stable process of type \textup{II}.

\mbox{} \textup{(ii)} In case $\alpha\in(0,1)$, we suppose that the assumptions of
either part \textup{(ii)} of Theorem
\ref{stable-discrete-1}, or of Theorem \ref{thm-coupl} hold. Then
the finite distributions of processes $\{Y_N(t), t\ge0\}$ converge,
as $N\to+\infty$,
to the respective distributions of a stable process of type \textup{I}.

\textup{(iii)} When $\alpha=1$ and the assumptions of part \textup{(iii)} of Theorem
\ref{stable-discrete-1} hold, the finite distributions of processes $\{
Y_N(t)-c_Nt, t\ge0\}$ converge,
as $N\to+\infty$,
to the respective distributions of a stable process of type \textup{III}.
Here $c_N:=\int_{|\Psi|\le N} \Psi \,d\pi$.
\end{theorem}

\section{An application: Superdiffusion of energy in a lattice dynamics}
\label{sec:appl-superd-energy}

In \cite{basollaspohn} it is proven that the Wigner distribution
associated with the energy of a system of interacting oscillators with
momentum and energy conserving noise converges, in an appropriate kinetic
limit, to the solution
$u(t,x,k)$ of the linear kinetic equation
%
%
\begin{equation}
\label{kinetic-eqt}
\cases{
\partial_t u(t,x,k)+\omega'(k)\,\partial_x u(t,x,k)=\mathcal{L}
u(t,x,k),\cr
u(0,x,k)=u_0(x,k),}
\end{equation}
where $(t,x,k)\in[0,+\infty)\times{\mathbb R}^d\times{\mathbb T}^d$
and the initial condition $u_0(\cdot,\cdot)$ is a function of class
$C^{1,0}({\mathbb R}^d\times{\mathbb T}^d)$.
Here ${\mathbb T}$ is the one-dimensional circle, understood as the
interval $[-1/2,1/2]$ with identified endpoints, and ${\mathbb T}^d$ is the
$d$-dimensional torus.
The function $\omega(k)$ is the dispersion relation of the lattice and it
is assumed that $\omega(-k)=\omega(k)$ and $\omega(k) \sim|k|$ for
$|k|\sim0$ (\textit{acoustic dispersion}). The scattering operator
$\mathcal{L}$, acting in (\ref{kinetic-eqt}) on variable $k$, is
usually an integral operator that is a generator of a certain jump process.

In the case of $d=1$, the scattering operator
is given by
%
%
\begin{equation}
\label{kinetic-eqt-1}
\mathcal{L}f(k)=\int_{{\mathbb T}}R(k,k')[f(k')-f(k)]\,dk'
\end{equation}
with the scattering kernel,
%
%
\begin{eqnarray}
\label{kernel}
R(k,k')&=& 
\tfrac43 [2\sin^2(2\pi k) \sin^2(\pi k')\nonumber\\[-8pt]\\[-8pt]
&&\hspace*{7.72pt}{} + 2\sin^2(2\pi k')
\sin^2(\pi k) - \sin^2(2\pi k) \sin^2(2 \pi k') ].\nonumber
\end{eqnarray}

We shall assume that the dispersion relation is given by a function
$\omega\dvtx\mathbb T\to[0,+\infty)$, that satisfies $\omega\in
C^1(\mathbb
T\setminus\{0\})$ and
%
%
\begin{equation}
\label{dispersion}
c_l|{\sin(\pi k)}|\le\omega(k)\le c_u|{\sin(\pi k)}|,\qquad k\in
\mathbb T,
\end{equation}
for some $0<c_l\le c_u<+\infty$ while
%
%
\begin{equation}
\label{dispersion-1}
\lim_{k\to\pm0} \omega'(k)=\pm c_\omega.
\end{equation}
In the case of a simple one-dimensional lattice, we have $\omega(k) =
c|{\sin(\pi k)}|$.

The total scattering cross section is given by
%
%
\begin{equation}
\label{eq:11}
R(k)=\int_{\mathbb T}R(k,k')\,dk' = \frac{4}{3} \sin^2(\pi k) \bigl(
1 +
2 \cos^2(\pi k) \bigr).
\end{equation}
We define $t(k) := R(k)^{-1}$ since these are the expected waiting
times of the scattering process.



Let $\{X_n, n\ge0\}$ be a Markov chain on $\mathbb T$ whose
transition probability equals
\[
P(k,dk'):= t(k) R(k,k') \,dk'.
\]
%
Suppose that $\{\tau_n, n\ge0\}$
is an i.i.d. sequence of random variables such that
$\tau_0$ is
exponentially distributed with intensity 1.
Let $t_n:=t(X_n)\tau_n$, $n\ge0$.
One can represent then the solution of (\ref{kinetic-eqt})
with the formula
%
%
\begin{equation}
\label{prob-rep}
u(t,x,k)=\mathbb E u_0(x(t),k(t)),
\end{equation}
where
\begin{eqnarray*}
x(t) &=& x+\int_0^t\omega'(k(s))\,ds,
\\
k(t)&=&X_n,\qquad t\in[t_n,t_{n+1}),
\end{eqnarray*}
and $k(0)=X_0=k$.
We shall be concerned in determining the weak limit of the
finite-dimensional distribution of the scaled process
$\{N^{-1/\alpha}x(Nt), t\ge0\}$, as $N\to+\infty$, for an appropriate
scaling exponent $\alpha$.

It is straightforward to verify that
%
%
\begin{equation}
\label{eq:12}
\pi(dk) = \frac{t^{-1}(k)}{\bar R} \,dk = \frac{R(k)}{\bar R} \,dk,
\end{equation}
where $\bar R:=\int_{{\mathbb T}} R(k)\,dk$ is a stationary and reversible
measure for the chain.
Then $P(k,dk')=p(k,k')\pi(dk')$ where
\[
p(k,k') = \bar R t(k) R(k,k') t(k')
\]
and after straightforward calculations, we obtain
%
%
\begin{eqnarray}
\label{eq:19}
p(k,k') &=&
6[\cos^2(\pi k)+\cos^2(\pi k')-
2\cos^2(\pi k)\cos^2(\pi k')]\nonumber\\
&&{}\times\bigl[ \bigl( 1 +
2 \cos^2(\pi k) \bigr)
\bigl( 1 +
2 \cos^2(\pi k') \bigr) \bigr]^{-1} \nonumber\\
&=&6\{[|{\cos(\pi k)}|-|{\cos(\pi k')}|]^2\\
&&\hspace*{7.93pt}{} +
2|{\cos(\pi k)\cos(\pi k')}|[1-|{\cos(\pi k)\cos(\pi k')}|]\}\nonumber\\
&&{}\times\bigl[ \bigl( 1 +
2 \cos^2(\pi k) \bigr)
\bigl( 1 +
2 \cos^2(\pi k') \bigr) \bigr]^{-1} .\nonumber
\end{eqnarray}
We apply Theorem \ref{thm-main-3} and probabilistic representation
(\ref{prob-rep}) to describe the asymptotic behavior for long times
and large spatial scales
of solutions of the kinetic equation (\ref{kinetic-eqt}). The result
is contained in the following.
\begin{theorem}
\label{kin-stable}
The finite-dimensional distributions of scaled processes
$ \{N^{-2/3}x(Nt), t\ge0 \}$ converge weakly to those
of a stable process of type \textup{II}. In addition, for any $t>0$, $x\in
{\mathbb R}
$, we have
%
%
\begin{equation}
\label{100801}
\lim_{N\to+\infty}
\int_{\mathbb T}|u(Nt,N^{2/3}x,k)-\bar u(t,x)|^2\,dk=0,
\end{equation}
where $u(t,x,k)$ satisfies
(\ref{kinetic-eqt}) with the initial condition $u_0(N^{-2/3}x,k)$,
such that $u_0$ is compactly supported, and $\bar u(t,x)$ is the
solution of
%
%
\begin{equation}
\label{fractional-heat}
\cases{
\partial_t \bar u(t,x)=-(-\partial^{2}_x)^{3/4}\bar u(t,x),\vspace*{2pt}\cr
\bar u(0,x)=\displaystyle\int_{\mathbb T} u_0(x,k)\,dk.}
\end{equation}
\end{theorem}
\begin{pf}
We start verifying the hypotheses of Theorem \ref{thm-main-3} by
finding the tails of
%
%
\begin{equation}
\label{082601}
\Psi(k) = \omega'(k) t(k)
\end{equation}
under measure $\pi$.
Since $\omega'(k)$ is both bounded and bounded away from zero, the
tails of $\Psi(k)$ under $\pi$ are the same as those of $t(k)$.
Note that
%
%
\begin{equation}
\label{eq:23}
\pi\bigl(k\dvtx t(k) \ge\lambda\bigr)=C_R \lambda^{-3/2} \bigl(1 +
O(1)\bigr) \qquad\mbox{for }\lambda\gg1,
\end{equation}
and some $C_R>0$. This verifies (\ref{eq:tails}) with $\alpha=3/2$.
Since the density of $\pi$ with respect to the Lebesgue measure is
even and $\Psi$ is odd, it has a null $\pi$-average.
\end{pf}

\subsection*{Verification of hypotheses of part \textup{(i)} of Theorem \protect
\ref{stable-discrete-1}}
Note
that we can decompose $P(k,dk')$ as in (\ref{transition}) with
$p(k,k')$ given by (\ref{eq:19}) and $Q(k,dk')\equiv0$. Since
$p(k,k')$ is bounded, Condition \ref{decomp} and (\ref{020411}) are
obviously satisfied.
Operator $P$ is a contraction on
$L^2(\pi)$, and by the Hilbert--Schmidt theorem (see, e.g., Theorem~4,
page 247 of \cite{lax}) is
symmetric and compact.
In consequence, 
its spectrum is contained in $[-1,1]$ and
is discrete, except for a possible accumulation point at~$0$.
\begin{lemma}
\label{lm:2}
Point $1$ is a simple eigenvalue of both $P$ and $P^2$.
\end{lemma}
\begin{pf}
Suppose
%
%
\begin{equation}\label{eq:10}
Pf=f.
\end{equation}
We claim that
$f$ is either everywhere positive, or everywhere negative.
Let $f^+,f^-$
be the positive and negative parts of $f$. Suppose also that $f^+$ is
nonzero on a set of positive $\pi$ measure.
Then
$f=f^+-f^-$
and $Pf=Pf^+-Pf^-$. Thus $f^+=(Pf)^+\le Pf^+$. Yet
\[
\int f^+\,d\pi\le\int Pf^+\,d\pi=\int f^+\,d\pi,
\]
thus $Pf^+=f^+$. Likewise, $Pf^-=f^-$. Since for each $k$ we have
$p(k,k')>0$, except
for a set of $k'$ of measure $\pi$ zero,
we conclude that $f^+>0$ $\pi$ a.e., hence $f^-\equiv0$.

Now we know that $P1=1$.
We claim that any other $f\not\equiv0$ that satisfies (\ref{eq:10})
belongs to $\operatorname{span}\{1\}$. Otherwise $f-c1$ for some $c$ would suffer
change of sign. But this contradicts our conclusion reached above so
the lemma holds for $P$.
The argument for $P^2$ is analogous.

As a corollary of the above lemma we conclude that
condition \ref{sg} holds.
Applying part (i) of Theorem \ref{thm-main-3}
to $N^{-2/3}\int_0^{Nt}\omega'(k(s))\,ds$, we conclude that its
finite-dimensional distributions converge in
law to an $\alpha$-stable L\'{e}vy process for $\alpha= 3/2$.

We use the above result to prove (\ref{100801}).
To
abbreviate the notation denote
$Y_N(t):=x+N^{-2/3}\int_0^{Nt}\omega'(k(s))\,ds$.
Using probabilistic representation for a solution of (\ref
{kinetic-eqt}), we
can write
%
%
\begin{eqnarray}
\label{extrac}\quad
u(Nt,N^{3/2}x,k)
&=&{\mathbb E}_k u_0 (Y_N(t),k(Nt) )\nonumber\\[-8pt]\\[-8pt]
&=&\sum_{\eta\in{\mathbb Z}}\int_{{\mathbb R}} \hat u_0(\xi,\eta
){\mathbb E}_k\exp
\{i\xi
Y_N(t)+i\eta k(Nt) \}\,d\xi.\nonumber
\end{eqnarray}
Here $\hat u_0(\xi,\eta)$ is the Fourier transform of $u(x,k)$,
and ${\mathbb E}_k$ is the expectation with respect to the path measure
corresponding to the momentum process $\{k(t), t\ge0\}$ that
satisfies $k(0)=k$. Since the dynamics of the momentum process are
reversible with respect to the normalized Lebesgue measure $m$ on the
torus and $0$ is a simple eigenvalue for the generator $\mathcal{L}$, we
have $\|P^tf\|_{L^2(m)}\to0$, as $t\to+\infty$, provided
$\int_{\mathbb T} fdk=0$. Suppose that $\{a_N, N\ge1\}$ is an
increasing sequence of positive numbers tending to infinity and such
that $a_NN^{-3/2}\to0$. A simple calculation shows that
for any $\xi,\eta\in\mathbb R$ and $e_\xi(x):=e^{ix\xi}$, we have
%
%
\begin{eqnarray}
\label{extra-b}\qquad
\bigl|{\mathbb E}_k[e_\xi( Y_N(t))e_\eta( k(Nt))]-{\mathbb
E}_k\bigl[e_\xi\bigl(
Y_N(t- ta_N/N)\bigr)e_\eta( k(Nt))\bigr] \bigr|\to0\nonumber\\[-8pt]\\[-8pt]
\eqntext{\mbox{as }N\to+\infty.}
\end{eqnarray}
Using Markov property we can write that the second term under the
absolute value in the formula above equals
\[
{\mathbb E}_k\bigl[e_\xi\bigl( Y_N(t- ta_N/N)\bigr)P^{a_Nt}e_\eta\bigl(k\bigl((N-a_N)t\bigr)\bigr)\bigr].
\]
Let $\tilde e_{\eta}(k):=e_\eta(k)-\bar e_\eta$, where $\bar e_{\eta
}:= \int_{\mathbb T} e_\eta(k)\,dk$.
By the Cauchy--Schwarz inequality, we obtain
%
%
\begin{eqnarray}
\label{extra-a}
&&\bigl|{\mathbb E}_k\bigl[e_\xi\bigl(
Y_N(t- ta_N/N)\bigr)P^{a_Nt}e_\eta\bigl(k\bigl((N-a_N)t\bigr)\bigr)\bigr]\nonumber\\
&&\hspace*{88.54pt}{} -{\mathbb E}_ke_\xi\bigl(
Y_N(t- ta_N/N)\bigr)
\bar e_\eta\bigr| \\
&&\qquad\le\bigl\{{\mathbb E}_k \bigl|P^{a_Nt}\tilde
e_\eta\bigl(k((N-a_N)t)\bigr) \bigr|^2 \bigr\}^{1/2}.\nonumber
\end{eqnarray}
The right-hand side of (\ref{extra-a}) tends to $0$ in the $L^2$ sense
with respect to $k\in\mathbb T$, as $N\to+\infty$.

From (\ref{extra-a}) we conclude that
%
%
\begin{eqnarray}
\label{extra-d}
&& \biggl|\sum_{\eta\in{\mathbb Z}}\int_{{\mathbb R}} \int
_{\mathbb T} \hat
u_0(\xi,\eta){\mathbb E}_k \bigl[e_\xi\bigl(
Y_N(t- ta_N/N)\bigr)\nonumber\\
&&\hspace*{96.18pt}{}\times P^{a_Nt}e_\eta\bigl(k\bigl((N-a_N)t\bigr)\bigr) \bigr]\,d\xi\,
dk
\\
&&
\hspace*{11.7pt}{}-\sum_{\eta\in{\mathbb Z}}\int_{{\mathbb R}} \int_{\mathbb
T} \hat
u_0(\xi,\eta){\mathbb E}_ke_\xi\bigl(
Y_N(t- ta_N/N)\bigr)\bar e_\eta \,d\xi \,d k \biggr|\to0
\nonumber
\end{eqnarray}
as $N\to+\infty$. Combining this with (\ref{extra-b}), we complete the
proof of the theorem.
\end{pf}

\subsection*{Verification of hypotheses of Theorem \protect\ref{thm-coupl}}

Here we show the convergence of $N^{-3/2}x(Nt)$ by using the coupling
approach of Section \ref{sec:coupling-approach}. Define the functions
\begin{eqnarray*}
q_0(k)&:=&\sin^2(2\pi k)=4 [\sin^2(\pi k) -\sin^4(\pi k) ],
\\
q_1(k)&:=& \tfrac{4}{3}\sin^4( \pi k),
\end{eqnarray*}
which are densities with respect to the Lebesgue measure in $\mathbb T$.
A simple computation shows that $R(k,k') =
2^{-4}[q_0(k)q_1(k')+q_1(k)q_0(k')]$, and therefore, $R(k) =
2^{-4}[q_0(k)+q_1(k)]$. The transition probability $P(k,dk')$ can be
written as
\[
P(k,dk') = \frac{q_1(k)}{q_0(k)+q_1(k)} q_0(k')\,dk' + \frac
{q_0(k)}{q_0(k)+q_1(k)} q_1(k')\,dk'.
\]

In particular, in the notation of Section \ref{sec:coupling-approach},
this model satisfies Condition \ref{doeblin} with $q(dk')
=q_0(k')\,dk'$, $\theta= q_1/(q_0+q_1)$ and $Q_1(k,dk') =q_1(k')\,dk'$.
Notice that the behavior around 0 of $\pi$
and $q$ is the same. Hence, $q(\Psi(k)\ge\lambda) \sim
c\lambda^{-3/2}$ for \mbox{$\lambda\gg1$}. We conclude, therefore, that the
function $\Psi(k)$, given by (\ref{082601}), satisfies~(\ref{eq:tails1}).
Observe furthermore that $Q_1$ does not depend on
$k$ and that
$Q_1(k',t(k)\ge\lambda)\sim c\lambda^{-5/2}$ for $\lambda\gg1$. Due
to this last observation, condition (\ref{tailsQ1}) is satisfied.

We are only left to check Condition \ref{regen}. But this one is also
simple, once we observe that the sequence $\{\delta_n, n\ge0\}$ is a
Markov chain with transition probabilities
\begin{eqnarray*}
P(\delta_{n+1}=1|\delta_n=0) &=& P(\delta_{n+1}=0|\delta_n = 1)
\\
&=& \int_{-1/2}^{1/2} \frac{q_0(k)q_1(k)\,dk}{q_0(k)+q_1(k)},\\
P(\delta_{n+1}=1|\delta_n=1) &=& \int_{-1/2}^{1/2} \frac
{q_1^2(k)\,dk}{q_0(k)+q_1(k)},
\\
P(\delta_{n+1}=0|\delta_n=0) &=& \int_{-1/2}^{1/2} \frac
{q_0^2(k)\,dk}{q_0(k)+q_1(k)}.
\end{eqnarray*}
We conclude that the regeneration time $\kappa_1$ satisfies
${\mathbb E}[\exp\{\gamma\kappa_1\}]<+\infty$ for $\gamma$ small
enough. Condition \ref{regen} is therefore a consequence of the fact
that the transition probability function $Q_1(k,dk')$ does not depend
on $k$; therefore, we can write
\[
{\mathbb P}[\kappa_1\ge n|K_0=k]
=\bigl(1-\theta(k)\bigr){\mathbb P}[\kappa_1\ge n-1].
\]

\section[Proof of Theorem 2.7 by coupling]{Proof of Theorem \protect\ref{thm-coupl} by coupling}
\label{sec:coupling-approach}

Because of its simplicity, we present first the proof of Theorem \ref
{thm-coupl} using a basic coupling argument.
Let us define
\begin{eqnarray*}
\varphi_i &=& \sum_{j=\kappa_i}^{\kappa_{i+1}-1} \Psi(X_j) ,
\\
\mathcal M(N) &=& \sup\{i \geq0; \kappa_i \leq N\}.
\end{eqnarray*}
Note that $\mathcal M(N)<+\infty$ a.s.
An alternative way of defining $\mathcal M(N)$ is demanding the inequality
$\kappa_{\mathcal M(N)} \leq N < \kappa_{\mathcal M(N)+1}$ to be satisfied.
Then, we have
%
%
\begin{equation}
\label{082501}
S_N = \sum_{i=0}^{\mathcal M(N)} \varphi_i + R_N,
\end{equation}
where
\[
R_N:=\sum_{j = \kappa_{\mathcal M(N)}+1}^N \Psi(X_j).
\]
In (\ref{082501}) we have decomposed $S_N$ into a random sum of i.i.d.
random variables, $\{\varphi_i, i\ge1\}$, and
two boundary terms: $\varphi_0$ and $R_N$. Notice also that
$\kappa_N-\kappa_1$ is a sum of i.i.d. random variables. Consequently,
the law of large numbers gives
%
%
\begin{equation}\label{eq:25}
\frac{\kappa_N}{N} \to\bar\kappa= \mathbb E(\kappa_2 - \kappa_1)
\quad\mbox{and}\quad
\frac{\mathcal M(N)}{N} \to\bar\kappa^{-1} = \bar\theta,
\end{equation}
a.s., as $N\to+\infty$.

Observe also that when $\alpha\in(1,2)$ and $\Psi$ is centered,
random variable $\varphi_1$ is also centered. Indeed, by the ergodic theorem
we have that a.s.
\[
0=\lim_{N\to+\infty}\frac{S_N}{N}=\lim_{N\to+\infty}\frac
{1}{\mathcal{M}(N)} \sum_{i=1}^{\kappa_{\mathcal{M}(N)}}\varphi
_i\times\frac{\mathcal{M}(N)}{N}=\mathbb E\varphi_1 \bar\theta,
\]
which proves that
%
%
\begin{equation}
\label{phi-1}
\mathbb E\varphi_1=0.
\end{equation}

The idea now is that under Conditions \ref{doeblin} and \ref{regen},
the random variable $\varphi_i$ is equal to $\Psi(X_{\kappa_i})$ plus
a term with lighter tails. Before stating this result, we need a
simple lemma.
\begin{lemma}
\label{l1}
Let $\zeta$ be a random variable such that
\[
\lim_{x \to\infty} x^\alpha{\mathbb P}(\zeta>x) =c^+,\qquad \lim
_{x \to
\infty} x^\alpha{\mathbb P}(\zeta<-x)=c^-.
\]

Let $\xi$ be such that $\lim_{x \to\infty} {\mathbb P}(|\xi
|>x)/{\mathbb P}
(|\zeta|>x)=0$. Then
%
%
\begin{equation}
\label{x-01}
\lim_{x \to\infty} x^\alpha{\mathbb P}(\zeta+\xi>x)=c^+,\qquad
\lim_{x
\to\infty} x^\alpha{\mathbb P}(\zeta+\xi<-x) =c^-.
\end{equation}
\end{lemma}
\begin{pf}
Without loss of generality, we just consider the first limit, the
second one follows considering $-\zeta$, $-\xi$. We will prove that
the $\liminf_{x\to\infty}$ of the previous expression is bigger than
$c_+$ and the $\limsup$ is smaller than $c_+$. We start with the upper
bound: for any $\epsilon>0$ there exists $x_0$ so that for $x\ge x_0$,
we have
\begin{eqnarray*}
x^\alpha{\mathbb P}(\zeta+\xi>x)
&\leq& x^\alpha{\mathbb P}\bigl(\zeta> (1-\epsilon)x\bigr) + x^\alpha
{\mathbb P}(\xi> \epsilon x)\\
&\leq& \frac{c_++\epsilon}{(1-2\epsilon)^\alpha} + \frac{{\mathbb
P}(|\xi|>
\epsilon x)}{{\mathbb P}(|\zeta|> \epsilon x)}\times\frac{c_++c_-
}{(\epsilon/2)^\alpha}.
\end{eqnarray*}
Now take above the upper limit, as $x \to+\infty$, to get
\[
\limsup_{x \to+\infty} x^\alpha{\mathbb P}(\zeta+\xi>x) \leq
\frac
{c_++\epsilon}{(1-2\epsilon)^\alpha}.
\]
Since $\epsilon$ is arbitrary, we have proved the upper bound. The
lower bound is very similar:
\begin{eqnarray*}
{\mathbb P}(\zeta+\xi>x)
&=& {\mathbb P}(\zeta+\xi>x, \xi> -\epsilon x) + {\mathbb P}(\zeta
+\xi>x, \xi
\leq-\epsilon x) \\
&\geq&{\mathbb P}\bigl(\zeta> (1+\epsilon)x , \xi> -\epsilon x\bigr) \\
&\geq&{\mathbb P}\bigl(\zeta> (1+\epsilon)x\bigr) - {\mathbb P}\bigl(\zeta>
(1+\epsilon) x ,
\xi\leq-\epsilon x\bigr) \\
&\geq&{\mathbb P}\bigl(\zeta> (1+\epsilon) x\bigr) -{\mathbb P}(\xi< -\epsilon x).
\end{eqnarray*}

Starting from this last expression, the same computations done for the
upper bound show that
\[
\liminf_{x \to+\infty} x^\alpha{\mathbb P}(\zeta+\xi>x) \geq
\frac
{c_+}{(1+2\epsilon)^\alpha}.
\]

Since $\epsilon>0$ is arbitrary, the lemma is proved for the first
expression in (\ref{x-01}). The second case can be done in the same fashion.
\end{pf}
\begin{lemma}
\label{lem:tailind}
Let $\Psi$ satisfy (\ref{eq:tails}) with constants
$c_*^+,c_*^-$ together with Conditions
\ref{doeblin} and \ref{regen}. Then the law of each
$\varphi_i$ satisfies
%
%
\begin{eqnarray}
\label{tails-b}
\lim_{\lambda\to+\infty} \lambda^\alpha{\mathbb P}(\varphi
_i>\lambda
) &=& c_*^+ \bar\theta^{-1},
\nonumber\\[-8pt]\\[-8pt]
\lim_{\lambda\to+\infty} \lambda^\alpha{\mathbb P}(\varphi_i <
-\lambda) &=& c_*^- \bar\theta^{-1}.\nonumber
\end{eqnarray}
%
\end{lemma}
\begin{pf}
The idea of the proof is simple. Random variable $\varphi_i$ is the
sum of a random variable with an $\alpha$-tail,
$\Psi(X_{\kappa_i})$, and a finite (but random) number of random
variables with lighter tails ($\Psi(X_{\kappa_i+1}), \ldots,
\Psi(X_{\kappa_{i+1}-1})$). By Condition \ref{regen}, the random
number can be efficiently controlled.
To simplify the notation, assume that $X_0$ is distributed according to
$q$, so the first block is also distributed like the other ones. Then
%
%
\begin{eqnarray}
\label{x-03}
{\mathbb P}\Biggl(\sum_{j=1}^{\kappa_1-1} \Psi(X_j)\geq t\Biggr)
&=& \sum_{n=1}^\infty{\mathbb P}\Biggl(\sum_{j=1}^{n-1} \Psi(X_j)
\geq t, \kappa_1=n\Biggr)
\nonumber\\[-8pt]\\[-8pt]
&\leq&\sum_{n=1}^\infty\sum_{j=1}^{n-1}
{\mathbb P}\bigl(\Psi(X_j) \geq t/(n-1),\kappa_1 = n\bigr).\nonumber
\end{eqnarray}
The probability under the sum appearing in the last expression can be
estimated by
%
%
\begin{eqnarray}
\label{x-02}
&&{\mathbb P}\bigl(\Psi(X_j) \geq t/(n-1), \delta_p=1,\forall p\le
j\bigr)\nonumber\\[-8pt]\\[-8pt]
&&\qquad={\mathbb E} \bigl[Q_1\bigl(X_{j-1},\Psi\ge t/(n-1)\bigr), \delta_p=1,
\forall
p\le j \bigr].
\nonumber
\end{eqnarray}
%
When $j\ge n/2$ we can use (\ref{tailsQ1}) to bound the expression on
the right-hand side of (\ref{x-02}) from above by
%
%
\begin{equation}
\label{x-06}
\frac{n^{\alpha}g(t/n)}{t^{\alpha}}{\mathbb P} [ \delta_p=1,
\forall
p\le j ]\le\frac{n^{\alpha}g(t/n)}{t^{\alpha}}{\mathbb
P} [ \kappa
_1\ge n/2 ].
\end{equation}
Here $g(x)$ is a bounded function that goes to 0, as $x \to
\infty$. On the other hand, when $j<n/2$, we rewrite
the probability appearing under the sum on the right-hand side of
(\ref{x-03})
using the Markov property. It equals
\[
{\mathbb E} \bigl[{\mathbb E} [\tilde\kappa_1=n-j |
\tilde
X_0=X_j ], \Psi( X_j) \geq t/(n-1), \delta_0=\cdots=\delta
_j=1 \bigr].
\]
Here $\{\tilde X_n, n\ge0\}$ is another copy of the Markov chain $\{
X_n, n\ge0\}$, and $\tilde\kappa_1$ is the respective stopping time
defined in correspondence to $\kappa_1$.
We can estimate this expression by
\begin{eqnarray*}
&& {\mathbb P}[\Psi(X_j) \geq t/(n-1)]\sup_{x}{\mathbb P}[\kappa
_1\ge
n-j|X_0=x]\\
&&\qquad
\le\frac{n^{\alpha}g(t/n)}{t^{\alpha}}\sup_{x}{\mathbb P}[\kappa
_1\ge
n/2|X_0=x].
\end{eqnarray*}
Summarizing, we have shown that the utmost left-hand side of
(\ref{x-03}) can be estimated by
\[
\sum_{n=1}^\infty\frac{g(t/n) n^{1+\alpha}}{t^\alpha}
\sup_x {\mathbb P}_x(\kappa_1 \geq n/2).
\]
We conclude that this expression
is $o(t^{-\alpha})$ by invoking Lebesgue dominated convergence theorem
and Condition
\ref{regen}. The negative tails are treated in the same
way. Therefore, $\varphi_0 -\Psi(X_0)$ has lighter tails than
$\Psi(X_0)$ itself. By Lemma \ref{l1}, the sum of a random variable
satisfying condition (\ref{eq:tails1}) and a random variable with lighter
tails also satisfies condition (\ref{eq:tails1}) for the same constants
$c^+ \bar\theta^{-1}$, $c^-\bar\theta^{-1}$.

At this point we are only left to recall the classical limit theorem
for i.i.d. random variables.
It follows that there exist
\begin{eqnarray*}
N^{-1/\alpha} S_N &=& \biggl(\frac{\mathcal M(N)}{N} \biggr)^{1/\alpha
} \frac
1{\mathcal M(N)^{1/\alpha}} \sum_{i=0}^{\mathcal M(N)} \varphi_i +
\frac
1{N^{1/\alpha}} \sum_{j = \kappa_{\mathcal M(N)}+1}^N \Psi(X_j)\\
&=& \biggl(\frac{\mathcal M(N)}{N\bar\theta} \biggr)^{1/\alpha}
\frac
1{\mathcal M(N)^{1/\alpha}} \sum_{i=0}^{\mathcal M(N)} \bar\theta
^{1/\alpha
}\varphi_i + \frac
1{N^{1/\alpha}} \sum_{j = \kappa_{\mathcal M(N)}+1}^N \Psi(X_j).
\end{eqnarray*}
Recall (\ref{eq:25}), and notice that by (\ref{tails-b}),
\begin{eqnarray*}
\lim_{\lambda\to+\infty} \lambda^\alpha
{\mathbb P}(\bar\theta^{1/\alpha}\varphi_i>\lambda)&=& c_*^+,
\\
\lim_{\lambda\to+\infty} \lambda^\alpha
{\mathbb P}(\bar\theta^{1/\alpha}\varphi_i < -\lambda)&=& c_*^-.
\end{eqnarray*}
%
Let $C_*:=(c_*^-+c_*^+)\bar\theta^{-1}$. By virtue of the stable limit
theorem for i.i.d. random variables (see, e.g., \cite{durrett}, Theorem
7.7, page 153), we know that for $c_N:=N\mathbb E[\varphi_1, |\varphi
_1|\le(C_*N)^{1/\alpha}]$ such that\vspace*{1pt} the laws of
$
N^{-1/\alpha} ( \sum_{i=0}^{N}
\bar\theta^{1/\alpha}\varphi_i-c_N )
$ converge to an
$\alpha$-stable law.
When $\alpha<1$ constants $c_N\sim cN^{1/\alpha}$ and they can be discarded.
Observe, however, that since $\mathbb E\varphi_1=0$, cf. (\ref
{phi-1}), for $\alpha\in(1,2)$, we have
\begin{eqnarray*}
c_N &=& -N\mathbb E[\varphi_1, |\varphi_1|> (C_*N)^{1/\alpha}]=N\int
_{(C_*N)^{1/\alpha}}^{+\infty}\bigl[{\mathbb P}[\varphi_1<-\lambda
]-{\mathbb P}
[\varphi_1>\lambda]\bigr]\,d\lambda\\
&=& \bar\theta^{-1}N\int_{(C_*N)^{1/\alpha}}^{+\infty
}(c_*^-c_*^+)\,\frac{d\lambda}{\lambda^{\alpha
}}=C\bigl(1+o(1)\bigr)N^{1/\alpha}
\end{eqnarray*}
for some constant $C$. The constants $c_N$ can again be discarded.
We conclude therefore that the laws of
\[
\mathcal K_N:=N^{-1/\alpha} \Biggl( \sum_{i=0}^{N}
\bar\theta^{1/\alpha}\varphi_i \Biggr)
\]
weakly converge to some
$\alpha$-stable law $\nu_*$.
Since $\mathcal{L}_N:=\bar\theta^{-1/\alpha}N^{-1}\mathcal{M}(N)$
converges a.s. to 1, the joint law of
$(\mathcal{K}_N,\mathcal{L}_N)$ converges to $\nu_*\otimes\delta
_1$, as $N\to+\infty$. According to the Skorochod representation
theorem there exists a probability space and random variables
$(\bar{\mathcal{K}}_N,\bar{\mathcal{L}}_N)$ such that
$(\bar{\mathcal{K}}_N,\bar{\mathcal{L}}_N)\stackrel{d}{=}(\mathcal
{K}_N,\mathcal{ L}_N)$ for each $N$ and $(\bar{\mathcal{K}}_N,\bar
{\mathcal{L}}_N)\to
(Y_*,1)$ a.s.
The above\vspace*{-2pt} in particular implies that $\bar{\mathcal{K}}_{N\bar
{\mathcal{L}}_N}$
converges a.s. to $Y_*$. Since $\bar{\mathcal{K}}_{N\bar{\mathcal
{L}}_N}\stackrel{d}{=}
\mathcal{K}_{N\mathcal{L}_N}$, we conclude the convergence of the
laws of
$\mathcal{K}_{N\mathcal{L}_N}$ to $\nu_*$.
\end{pf}

\section[Proof of Theorem 2.4 by
martingale approximation]{Proof of Theorem \protect\ref{stable-discrete-1} by
martingale approximation}
\label{sec5}
Below we formulate a
stable limits law that shall be crucial during the course of the proof
of the theorems.

Suppose that $\{Z_n, n\ge1\}$ is a stationary sequence that is
adapted with respect to the filtration $\{\mathcal{G}_n, n\ge0\}$
and such that
for any $f$ bounded and measurable, the sequence $ \{\mathbb
E[f(Z_n)|\mathcal{G}_{n-1}], n\ge1 \}$ is also stationary. We
assume furthermore
that there exist $\alpha\in(0,2)$ and $c_*^+,c_*^-\ge0$ such that
$c_*^++c_*^->0$
and
%
%
\begin{eqnarray}
\label{eq:a2}
\mathbb P[Z_1>\lambda] &=&
\lambda^{-\alpha}\bigl(c_*^++o(1)\bigr),\nonumber\\[-8pt]\\[-8pt]
\mathbb P[Z_1<-\lambda] &=& \lambda^{-\alpha}\bigl(c_*^-+o(1)\bigr)\qquad
\mbox{as } \lambda\to+\infty.\nonumber
\end{eqnarray}
In addition, for any $g\in C_0^\infty(\mathbb R\setminus\{0\})$, we have
%
%
\begin{equation}
\label{eq:a3}
\lim_{N\to+\infty}\mathbb E \Biggl|\sum_{n=1}^{[Nt]}\mathbb E
\biggl[g \biggl(
\frac{Z_{n}}{N^{1/\alpha}} \biggr) \bigg|\mathcal
{G}_{n-1} \biggr]-\alpha t\int_{\mathbb R}g(\lambda)\frac
{c_*(\lambda) \,d\lambda
}{|\lambda|^{1+\alpha}} \Biggr|=0
\end{equation}
and
%
%
\begin{equation}
\label{eq:a4}
\lim_{N\to+\infty}N\mathbb E \biggl\{\mathbb E \biggl[g \biggl(
\frac{Z_{1}}{N^{1/\alpha}} \biggr) \bigg|\mathcal{G}_0
\biggr] \biggr\}^2=0.
\end{equation}
Here $c_*(\cdot)$ appearing in (\ref{eq:a3}) is given by (\ref{c-levy}).
Let $M_N:=\sum_{n=1}^{N}Z_n$, $N\ge1$ and $M_0:=0$.

When $\alpha=1$ we shall also consider
an array $\{Z_n^{(N)}\dvtx n\ge1\}$, $N\ge1$ of stationary sequences
adapted with respect to the filtration $\{\mathcal{G}_n\dvtx n\ge0\}$.
Assume furthermore that
for each $N\ge1$ and any $f$ bounded and measurable sequence, $ \{
\mathbb E[f(Z_n^{(N)})|\mathcal{G}_{n-1}]\dvtx n\ge1 \}$ is stationary.
We suppose that there exist nonnegative $c_*^+,c_*^-$ such that
$c_*^++c_*^->0$ and
%
%
\begin{equation}
\label{eq:a2b}
\lim_{\lambda\to+\infty}\sup_{N\ge1}\bigl[\bigl|\lambda\mathbb
P\bigl[Z_1^{(N)}>\lambda
\bigr]-c_*^+\bigr|+\bigl|\lambda\mathbb P\bigl[Z_1^{(N)}<-\lambda\bigr]-c_*^-\bigr|\bigr]=0.
\end{equation}
Let
\[
\tilde M_N:=\sum_{n=1}^{N} \bigl\{Z_n^{(N)}-{\mathbb E}
\bigl[Z_n^{(N)}1\bigl[\bigl|Z_n^{(N)}\bigr|\le N\bigr]\big|\mathcal{G}_{n-1}\bigr] \bigr\},\qquad
N\ge1,
\]
and $\tilde M_0:=0$.
The following result has been shown in Section 4 of \cite{durrettresnick}.
\begin{theorem}
\label{thm:2}
\textup{(i)}
Suppose that $\alpha\in(1,2)$, conditions (\ref{eq:a2})--(\ref{eq:a4})
hold, and
%
%
\begin{equation}
\label{eq:a1}
\mathbb E[Z_{n}|\mathcal{G}_{n-1}]=0 \qquad\mbox{for }n\ge1.
\end{equation}
Then
$N^{-1/\alpha}M_{[N\cdot]}\Rightarrow Z(\cdot)$, as $N\to+\infty$,
weakly in $D[0,+\infty)$,
where $\{Z(t),
t\ge0\}$ is an $\alpha$-stable
process of type \textup{II}.

\mbox{}\hspace*{2.8pt}\textup{(ii)} Suppose that $\alpha\in(0,1)$ and conditions (\ref{eq:a2})--(\ref
{eq:a4}) hold. Then
$N^{-1/\alpha}\times\break M_{[N\cdot]}\Rightarrow Z(\cdot)$, as $N\to+\infty$,
weakly in $D[0,+\infty)$,
where $\{Z(t),
t\ge0\}$ is an $\alpha$-stable
process of type \textup{I}.

\textup{(iii)} For\vspace*{1pt} $\alpha=1$, assume (\ref{eq:a3}) and (\ref{eq:a4}) with
$Z_n^{(N)}$ replacing $Z_n$ and (\ref{eq:a2b}).
Then $N^{-1}\tilde M_{[N\cdot]}\Rightarrow Z(\cdot)$, as $N\to
+\infty$, weakly in $D[0,+\infty)$
to a L\'{e}vy process $\{Z(t),
t\ge0\}$ of type \textup{III}.
\end{theorem}

\subsection*{Proof of part \textup{(i)} of Theorem \protect\ref{stable-discrete-1}}

Let $\chi\in L^\beta(\pi)$, $\beta\in(1,\alpha)$ be the unique, zero-mean
solution of the equation
%
%
\begin{equation}\label{eq:14}
\chi-P\chi=\Psi.
\end{equation}
Since $\Psi\in L^\beta(\pi)$ for $\beta\in(0,\alpha)$ is of zero
mean, the
solution to (\ref{eq:14}) exists in $L^\beta(\pi)$ and is given by
$\chi=\sum_{n\ge0}P^{n}\Psi$. This follows from the fact that
$\|P^n\Psi\|_{L^\beta}\le a^{(2/\beta-1)n}\|\Psi\|_{L^\beta}$,
$n\ge0$ [see
(\ref{eq:9})], so the series defining $\chi$ geometrically
converges. Uniqueness is a consequence of (\ref{eq:9}). Indeed, if
$\chi_1$ was another zero-mean solution to (\ref{eq:14}), then
\[
\|\chi-\chi_1\|_{L^\beta}=\|P(\chi-\chi_1)\|_{L^\beta}\le
a^{1-|2/\beta-1|}\|\chi-\chi_1\|_{L^\beta},
\]
which clearly is possible only when $\chi-\chi_1=0$ (recall that $a<1$).
Note also that
from (\ref{020411}) it follows that in fact
$P\chi=(I-P)^{-1}(P\Psi)\in L^{\alpha'}(\pi)$. Thus in particular,
\begin{equation}
\label{decay}
\pi(|P\chi|>\lambda)\le\frac{\|P\chi\|_{L^{\alpha'}(\pi
)}^{\alpha'}}{\lambda
^{\alpha'}},
\end{equation}
and consequently $\chi$ satisfies the same tail condition as $\Psi$
[cf. (\ref{eq:tails})].

Then by using (\ref{eq:14}), we can write
%
%
\begin{equation}
\label{eq:5}
S_N = \sum_{n=1}^N \Psi(X_n) = \sum_{n=1}^N Z_n + P\chi(X_0) -
P\chi(X_N)
\end{equation}
with $Z_n = \chi(X_{n})-P\chi(X_{n-1})$.

In what follows, we denote
by $C_0^\infty(\mathbb R\setminus\{0\})$ the space of all $C^\infty$
functions
that are compactly supported in $\mathbb R\setminus\{0\}$.
According to part (i) of Theorem \ref{thm:2}, we only need to
demonstrate the following.
\begin{prop}
\label{prop:2}
For any $g\in C_0^\infty(\mathbb R\setminus\{0\})$, equalities (\ref
{eq:a3}) and (\ref{eq:a4}) hold.
\end{prop}

More explicitly, we have
\[
E \biggl[g \biggl(\frac{Z_{n}}{N^{1/\alpha}} \biggr) \bigg|\mathcal{
G}_{n-1} \biggr]
= \int g\bigl(N^{-1/\alpha} [\chi(y) - P\chi(X_{n-1})]\bigr) P(X_{n-1}, dy)
\]
and using the stationarity of $\pi$, we can bound the left-hand side of
(\ref{eq:a3}) by
%
%
\begin{eqnarray}
\label{082701}\qquad
&&\mathbb E \Biggl|\sum_{n=1}^{[Nt]} \int\biggl[g \bigl(
N^{-1/\alpha}[\chi(y) - P\chi(X_{n-1})] \bigr) -
g \biggl(\frac{\chi(y)}{N^{1/\alpha}} \biggr) \biggr] P(X_{n-1}, dy)
\Biggr| \nonumber\\
&&\qquad{} + \mathbb E \Biggl|\sum_{n=1}^{[Nt]} \int
g \biggl(\frac{\chi(y)}{N^{1/\alpha}} \biggr)P(X_{n-1}, dy) -
[Nt] \int g \biggl(\frac{\chi(y)}{N^{1/\alpha}} \biggr)\pi(dy)
\Biggr|
\\
&&\qquad{} + \biggl|[Nt]\int g \biggl(\frac{\chi(y)}{N^{1/\alpha}} \biggr)\pi(dy)
-\alpha t \int_{\mathbb
R}g(\lambda)\frac{c_*(\lambda)\,d\lambda}{|\lambda|^{1+\alpha
}} \biggr|,\nonumber
\end{eqnarray}
so (\ref{eq:a3}) is a consequence of the following three lemmas, each
taking care of the respective term of (\ref{082701}):
\begin{lemma}
\label{lm:apro}
%
%
\begin{eqnarray}\label{eq:6}
&&\lim_{N\to\infty} N\iint\bigl|g\bigl(N^{-1/\alpha} [\chi(y) - P\chi
(x)]\bigr) \nonumber\\[-8pt]\\[-8pt]
&&\hspace*{94.61pt}{}
- g (N^{-1/\alpha} \chi(y) ) \bigr| P(x, dy) \pi(dx) = 0.\nonumber
\end{eqnarray}
\end{lemma}
\begin{lemma}
\label{lm:4}
%
%
\begin{equation}\qquad
\label{eq:18}
\mathbb E \Biggl|\sum_{n=1}^{N}\mathbb E [g (N^{-1/\alpha}
\chi(X_{n}) ) |\mathcal{G}_{n-1} ]
- N\int g (N^{-1/\alpha} \chi(y) )\pi(dy)
\Biggr|=0.
\end{equation}
\end{lemma}
\begin{lemma}
\label{lm:3}
%
%
\begin{equation}
\label{eq:21}
\lim_{N\to\infty} \biggl|N\int
g (N^{-1/\alpha}\chi(y) ) \pi(dy)
-\alpha\int_{\mathbb R}g(\lambda)\frac{c_*(\lambda)\,d\lambda
}{|\lambda|^{1+\alpha
}} \biggr|=0,
\end{equation}
where
$c_*(\cdot)$ is given by (\ref{c-levy}).
\end{lemma}

\subsection*{Proof of Lemma \protect\ref{lm:apro}}
Suppose that $\operatorname{supp}g\subset[-M,M]\setminus[-m,m]$ for some
$0<m<M<+\infty$ and $\theta>0$.
Denote
\[
A_{N,\theta} = \{ (x,y)\dvtx | \chi(y) - \theta P\chi(x)| >
N^{1/\alpha} m \}.
\]
The left-hand side of (\ref{eq:6}) can be bounded from above by
\begin{eqnarray*}
&&N^{1-1/\alpha}\int_0^1 d\theta\iint\bigl|g'\bigl(N^{-1/\alpha} [\chi(y)
- \theta P\chi(x)]\bigr) P\chi(x) \bigr| P(x, dy) \pi(dx)\\
&&\qquad\le C N^{1-1/\alpha}\int_0^1 d\theta\iint_{A_{N,\theta}}
| P\chi(x) | P(x, dy) \pi(dx)\\
&&\qquad\le C N^{1-1/\alpha}\int_0^1 d\theta\biggl(\iint_{A_{N,\theta}}
P(x, dy) \pi(dx) \biggr)^{1- 1/\alpha'} \|P\chi\|_{L^{\alpha'}}.
\end{eqnarray*}
From the tail behavior of $\chi$ and of $P\chi$, [see (\ref{decay})
and the remark below that estimates], it is easy to see that
for any $\theta\in(0,1)$,
\begin{eqnarray*}
\iint_{A_{N,\theta}} P(x, dy) \pi(dx) &\le& {\mathbb P}[|\chi
(X_1)|\ge
(mN^{1/\alpha})/2] \\
&&{} +
{\mathbb P}[\theta|P\chi(X_0)|\ge(mN^{1/\alpha})/2] \\
&\le& C[(Nm^{\alpha})^{-1}+(Nm^{\alpha})^{-\alpha'/\alpha}]\\
&=&\frac
{C}{N}\bigl(1+o(1)\bigr)
\end{eqnarray*}
as $N\gg1$.
Since $\alpha' >\alpha$ we obtain (\ref{eq:6}).

\subsection*{Proof of Lemma \protect\ref{lm:4}}
To simplify the notation we assume that $\operatorname{supp}g\subset[m,M]$ for
$0<m<M<+\infty$.
Denote $B_{N,\lambda}= \{ y\dvtx \chi(y) \ge N^{1/\alpha}\lambda\}$.
We can rewrite the left-hand side of (\ref{eq:18}) as
%
%
\begin{equation}
\label{082801}
\mathbb E \Biggl| \int_0^\infty g'(\lambda) \sum_{n=1}^N
G_N(X_{n-1},\lambda) \,d\lambda\Biggr|,
\end{equation}
where
\[
G_N(x,\lambda) = P(x, B_{N,\lambda}) - \pi(B_{N,\lambda}) .
\]
Notice that $ \int G_N(y,\lambda) \pi(dy) = 0$ and
%
%
\begin{eqnarray}\label{Gbounds}
\int G_N^2(y, \lambda) \pi(dy)
&=&
\int P^2(y, B_{N,\lambda}) \pi(dy) - \pi^2(B_{N,\lambda})
\nonumber\\
&\le& 2 \int\biggl( \int_{B_{N,\lambda}} p(y,x)
\pi(dx) \biggr)^2 \pi(dy) \\
&&{} + 2 \int
Q^2(y, B_{N,\lambda}) \pi(dy) - \pi^2(B_{N,\lambda}).\nonumber
\end{eqnarray}
To estimate the first term on the utmost right-hand side,
we use the Cauchy--Schwarz inequality, while for the second one we
apply condition (\ref{new}). For $\lambda\ge m$, we can bound the
expression on the right-hand side of (\ref{Gbounds}) by
%
%
\begin{eqnarray}\label{Gbounds1}
&& \frac12 \pi(B_{N,m}) \int\!\!\int_{B_{N,m}} p^2(x,y)
\pi(dx) \pi(dy) + C \pi^2(B_{N,m})\nonumber\\[-8pt]\\[-8pt]
&&\qquad\le\frac{1}{N} o(1)\qquad \mbox{as $N\to\infty$,}\nonumber
\end{eqnarray}
by virtue of (\ref{eq:4}) and the remark after (\ref{decay}). Thus we
have shown that
%
%
\begin{equation}
\label{Gbounds2}
N\sup_{\lambda\ge m} \int G_N^2(y, \lambda) \pi(dy) \to0
\end{equation}
as $N\to\infty$.
We will show now that (\ref{Gbounds2}) and the spectral gap together imply
that
%
%
\begin{equation}
\label{G-la}
\sup_{\lambda\ge m}\mathbb E\Biggl| \sum_{n=1}^N G_N(X_{n-1},\lambda)\Biggr|^2
\to0
\end{equation}
as
$N\to\infty$.
Since $\operatorname{supp}g'\subset[m,M]$ expression in (\ref{082801}) can be
then estimated by
\[
\sup_{\lambda\ge m} \mathbb E \Biggl| \sum_{n=1}^N
G_N(X_{n-1},\lambda) \Biggr|\times\int_0^\infty|g'(\lambda
)|\,d\lambda\to0
\]
as $N\to+\infty$ and the conclusion of the lemma follows.

To prove (\ref{G-la}) let $u_N(\cdot,\lambda) = (I-P)^{-1} G_N(\cdot
,\lambda)$. By the spectral
gap condition~(\ref{eq:9}), we have
%
%
\begin{equation}
\label{082802}
\int u_N^2(y,\lambda) \pi(dy) \le\frac1{1-a^2}\int G_N^2(y,
\lambda) \pi(dy).
\end{equation}
We can then rewrite
\[
\sum_{n=1}^N G_N(X_{n-1},\lambda) = u_N(X_0) - u_N(X_N) + \sum_{n=1}^{N-1}U_n
,
\]
where $U_n=u_N(X_n) - Pu_N(X_{n-1})$, $n\ge1$ is a stationary sequence
of martingale differences with respect to the natural filtration
corresponding to $\{X_n, n\ge0\}$. Consequently,
\[
\mathbb E\Biggl| \sum_{n=1}^N G_N(X_{n-1},\lambda)\Biggr|^2 \le C N \int
u_N^2(y,\lambda)
\pi(dy) \to0
\]
and (\ref{G-la}) follows from (\ref{Gbounds2}) and (\ref{082802}).

\subsection*{Proof of Lemma \protect\ref{lm:3}}
To avoid long notation, we again assume that $\operatorname{supp}g\subset[m,M]$
for $0<m<M<+\infty$. The proof
in the case of $ g\subset[-M,-m]$ is virtually the same. Note that
\begin{eqnarray*}
&&
N\int g \biggl(
\frac{\chi(y)}{N^{1/\alpha}} \biggr)\pi(dy)\\
&&\qquad
=N\int\!\!\int_0^{+\infty} N^{-1/\alpha}g' \biggl(
\frac{\lambda}{N^{1/\alpha}} \biggr)1_{[0,\chi(y)]}(\lambda)\pi
(dy)\,d\lambda
\\
&&\qquad
= N\int_0^{+\infty}N^{-1/\alpha} g' \biggl(
\frac{\lambda}{N^{1/\alpha}} \biggr)\pi(\chi> \lambda)\,d\lambda
\\
&&\qquad
=N \int_0^{+\infty} g' (\lambda)\pi( \chi\ge
N^{1/\alpha}\lambda
)\,d\lambda.
\end{eqnarray*}
Thanks to (\ref{eq:tails}) the last expression tends, however, as
$N\to+\infty$, to
\[
\int_0^{+\infty}
g' (\lambda)\frac{c_*^+\,d\lambda}{\lambda^{\alpha
}}=\alpha
\int_{\mathbb R} g (\lambda)\frac{c_*(\lambda)\,d\lambda
}{|\lambda
|^{\alpha+1}}.
\]

\subsection*{Proof of Proposition \protect\ref{prop:2}}
We have already shown (\ref{eq:a3}), so only (\ref{eq:a4}) requires
a proof. To
simplify the notation we assume $Q\equiv0$.
Suppose that $\operatorname{supp}g\subset[m,M]$ for
some $0<m<M$. We can write
%
%
\begin{eqnarray}\label{021901}
&&
\mathbb E \biggl[g \biggl(
\frac{Z_{1}}{N^{1/\alpha}} \biggr) \bigg|\mathcal{G}_0 \biggr]
\nonumber\\
&&\qquad= \int g (N^{-1/\alpha}\Psi(y) )p(X_0,y)\pi(dy)
\nonumber\\[-8pt]\\[-8pt]
&&\qquad\quad{}+ N^{-1/\alpha}\int\!\!\int_0^1h(X_0,y)
g' \bigl(N^{-1/\alpha}\bigl(\Psi(y)+\theta h(X_0,y)\bigr) \bigr)\nonumber\\
&&\hspace*{98.29pt}{}\times p(X_0,y)\pi
(dy)\,d\theta,\nonumber
\end{eqnarray}
where $h(x,y):=P\chi(y)-P\chi(x)$.
Denote by $K_1$ and $K_2$ the first and the second terms appearing on
the right-hand side above.
By Cauchy--Schwarz inequality,
%
%
\begin{eqnarray}
\label{011901}
{\mathbb E}K_2^2
&\le&
\frac{\|g'\|_{\infty}^2}{N^{2/\alpha}}
\biggl[{\mathbb E}
\biggl(\int|P\chi(y)|p(X_0,y)\pi(dy) \biggr)^2
\nonumber\\
&&\hspace*{34.62pt}{}
+{\mathbb E} \biggl(\int|P\chi(X_0)|p(X_0,y)\pi(dy) \biggr)^2
\biggr]
\\
&\le&\frac{2\|g'\|_{\infty}^2\|P\chi\|_{L^2(\pi)}^2}{N^{2/\alpha
}}.\nonumber
\end{eqnarray}
Hence
$\lim_{N\to+\infty}N{\mathbb E}K_2^2=0$.

On the other hand,
%
%
\begin{equation}
\label{eq:31}
K_1\le\|g\|_\infty\int p(X_0,y)1[|\Psi(y)|>mN^{1/\alpha}/2]\pi(dy),
\end{equation}
and in consequence, by Jensen's inequality,
\begin{eqnarray*}
({\mathbb E}K_1^2 )^{1/2}
&\le&\|g\|_\infty\int({\mathbb E}
p^2(X_0,y))^{1/2}1[|\Psi(y)|>aN^{1/\alpha}/2]\pi(dy)
\\
&\le&\|g\|_\infty\biggl[\int\!\!\int p^2(x,y)1[|\Psi(y)|>aN^{1/\alpha
}/2]\pi(dx)\pi(dy) \biggr]^{1/2}\nonumber\\
&&{}\times\pi^{1/2}[|\Psi|>aN^{1/\alpha}/2].
\end{eqnarray*}
Thus we have shown that
%
%
\begin{eqnarray}
\label{051901}\quad
N{\mathbb E}K_1^2
&\le& N\pi[|\Psi|>aN^{1/\alpha}/2]\nonumber\\[-8pt]\\[-8pt]
&&{}
\times\int\!\!\int p^2(x,y)1[|\Psi(y)|>aN^{1/\alpha}/2]\pi(dx)\pi
(dy)\to
0.\nonumber
\end{eqnarray}
Condition (\ref{eq:a4}) is then a consequence of (\ref{011901}) and
(\ref{051901}).

\subsection*{Proof of part \textup{(ii)} of Theorem \protect\ref
{stable-discrete-1}}

The proof of this part relies on
part (ii) of Theorem \ref{thm:2}.
The following
analogue of Proposition \ref{prop:2} can be established.
\begin{prop}
\label{prop:2b}
Suppose that $\alpha\in(0,1)$. Then
for any $g\in C_0^\infty(\mathbb R\setminus\{0\})$,
%
%
\begin{equation}
\label{eq:16bb}
\lim_{N\to+\infty}\mathbb E \Biggl|\sum_{n=1}^{[Nt]}\mathbb
E \biggl[g \biggl(\frac{\Psi(X_n)}{N^{1/\alpha}}
\biggr) \bigg|\mathcal{G}_{n-1} \biggr]- t\int_{\mathbb
R}g(\lambda)\frac{C_*(\lambda) \,d\lambda}{|\lambda|^{1+\alpha
}} \Biggr|=0
\end{equation}
and
%
%
\begin{equation}
\label{eq:16bx}
\lim_{N\to+\infty}N\mathbb E \biggl\{\mathbb E \biggl[g \biggl(
\frac{\Psi(X_{1})}{N^{1/\alpha}} \biggr) \bigg|\mathcal
{G}_0 \biggr] \biggr\}^2=0.
\end{equation}
\end{prop}
\begin{pf}
The proof of this proposition is a simplified version of the argument
used in the proof of Proposition \ref{prop:2}. The expression in
(\ref{eq:16bb}) can be estimated by
%
%
\begin{eqnarray}
\label{082701b}
&&\mathbb E \Biggl|\sum_{n=1}^{[Nt]} \int
g \biggl(\frac{\Psi(y)}{N^{1/\alpha}} \biggr)P(X_{n-1}, dy) -
[Nt] \int g \biggl(\frac{\Psi(y)}{N^{1/\alpha}} \biggr)\pi(dy)\Biggr|
\nonumber\\[-8pt]\\[-8pt]
&&\qquad{} + \biggl|[Nt]\int g \biggl(\frac{\Psi(y)}{N^{1/\alpha}} \biggr)\pi(dy)
-\alpha t \int_{\mathbb
R}g(\lambda)\frac{c_*(\lambda)\,d\lambda}{|\lambda|^{1+\alpha
}} \biggr|.\nonumber
\end{eqnarray}
The proof that both the terms of the sum above vanish goes along the
lines of the
proofs of Lemmas \ref{lm:4} and \ref{lm:3}. We can repeat word by
word the argument used there, replacing this time $\chi$ by $\Psi$.
As for the proof of (\ref{eq:16bx})
it is identical with the respective part of the proof
of (\ref{eq:a4}) (the one concerning term $K_1$).
\end{pf}
%

\subsection*{Proof of part \textup{(iii)} of Theorem \protect\ref
{stable-discrete-1}}

Recall that $\Psi_N:=\Psi1[|\Psi|\le N]$.
Let $\chi_N$ be the unique, zero mean
solution of the equation
%
%
\begin{equation}\label{eq:14b}
\chi_N-P\chi_N=\Psi_N-c_N.
\end{equation}
We can then write,
%
%
\begin{equation}\hspace*{25pt}
\label{eq:5b}
S_N-Nc_N = \sum_{n=1}^N \bigl(\Psi(X_n)-c_N\bigr) = \sum_{n=1}^N Z_n^{(N)} +
P\chi_N(X_0) -
P\chi_N(X_N)
\end{equation}
with
\[
Z_n^{(N)} = \chi_N(X_{n})-P\chi_N(X_{n-1})+\Psi(X_n)1[|\Psi(X_n)|>N].
\]
We verify first assumptions (\ref{eq:a3}), (\ref{eq:a4}) and
(\ref{eq:a2b}).

Condition (\ref{eq:a2b}) is an obvious consequence of
the fact that
%
%
\begin{equation}
\label{031023}
Z_n^{(N)} = P\chi_N(X_{n})-P\chi_N(X_{n-1})+\Psi(X_n)-c_N
\end{equation}
and assumption (\ref{020411b}).
To verify the remaining hypotheses, suppose that $\operatorname{supp}g\subset(m,M)$
and $m<1<M$.
Let us fix $\delta>0$, to be further chosen later on, such that
$m<1-\delta<1+\delta<M$. We can then write $g=g_1+g_2+g_3$ where
each $g_i\in C^\infty({\mathbb R})$, $\|g_i\|_{\infty}\le\|g\|
_\infty$,
and the supports of $g_1,g_2,g_3$ are correspondingly contained in
$(m,1-\delta), (1-\delta,1+\delta), (1+\delta,M)$.
We prove (\ref{eq:a3}) and (\ref{eq:a4}) for
each of the function $g_i$s separately.
Note that
\[
{\mathbb E} \biggl[g_i \biggl(\frac{Z_{n}^{(N)}}{N} \biggr)
\bigg|\mathcal{
G}_{n-1} \biggr]
= \int g_i \bigl(w^{(N)}(X_{n-1},y) \bigr) P(X_{n-1}, dy),
\]
where
\[
w^{(N)}(x,y):=N^{-1}\Psi(y)1[|\Psi(y)|>N]+N^{-1} [\chi_N(y) - P\chi_N(x)].
\]
For $i=1$ and $i=3$, we essentially estimate in the same way as in
parts (i) and (ii) of the proof of the theorem, respectively. We shall
only consider here the case $i=2$.

Note that then $w^{(N)}(x,y)=w_1^{(N)}(x,y)$ where
%
%
\begin{equation}
\label{102301}
w_\theta^{(N)}(x,y)=N^{-1}[\Psi(y)-c_N]+N^{-1}\theta R_N(x,y)
\end{equation}
with $R_N(x,y):=P\chi_N(y) - P\chi_N(x)$. However,
\begin{eqnarray*}
g_2 \bigl(w^{(N)}(X_{n-1},y) \bigr)&=&
g_2 \bigl(N^{-1}\bigl(\Psi(y)-c_N\bigr) \bigr)\\
&&{} + N^{-1} R_N(X_{n-1},y)\int_0^1
g_2' \bigl(w_\theta^{(N)}(X_{n-1},y) \bigr)\,d\theta
\end{eqnarray*}
and
\begin{eqnarray*}
&&
{\mathbb E} \Biggl|\sum_{n=1}^N\int g_2 \bigl(w^{(N)}(X_{n-1},y) \bigr)
P(X_{n-1}, dy) \Biggr|\\
&&\qquad
\le
{\mathbb E} \Biggl|\sum_{n=1}^N\int g_2 \bigl(N^{-1}\bigl(\Psi
(y)-c_N\bigr) \bigr)
P(X_{n-1}, dy) \Biggr|\\
&&\qquad\quad{}
+\int\!\!\int\!\!\int_0^1 \bigl|g_2' \bigl(w_\theta^{(N)}(x,y) \bigr)R_N(x,y)\bigr|
P(x, dy)\pi(dy)\,d\theta.
\end{eqnarray*}
Denote the first and the second term on the right-hand side by
$J_1^{(N)}$ and $J_2^{(N)}$, respectively. Term $J_1^{(N)}$ can be now
estimated as in the proof of part (ii) of the theorem. We conclude then,
using the arguments contained in the proofs of Lemmas~\ref{lm:4} and \ref{lm:3} that
\[
\limsup_{N\to+\infty}J_1^{(N)}\le\|g\|_\infty\int_{1-\delta
}^{1+\delta}\frac{d\lambda}{\lambda^2}.
\]
On the other hand, to estimate $\lim_{N\to+\infty}J_2^{(N)}=0$, since
$g_2' (w_\theta^{(N)}(x,y) )\to0$ in measure $P(x, dy)\pi
(dy)\,d\theta$
and the passage to the limit under he integral can be substantiated
thanks to
(\ref{020411b}).

Choosing now sufficiently small $\delta>0$ we can argue that the
calculation of the limit can be reduced to the cases considered for
$g_1$ and $g_3$ and that condition (\ref{eq:a3}) can be established
for $Z_n^{(N)}$.
The proof of (\ref{eq:a4}) can be repeated from the argument for part
(i) of the theorem.

Finally, we show that
%
%
\begin{equation}
\label{021023}
\lim_{N\to+\infty}\frac{1}{N}{\mathbb E} \Biggl|\sum
_{n=1}^{N}{\mathbb E}
\bigl[Z_n^{(N)}1\bigl[\bigl|Z_n^{(N)}\bigr|\le N\bigr]\big|\mathcal{G}_{n-1}\bigr] \Biggr|=0.
\end{equation}
Denote the expression under the limit by $L^{(N)}$. Let $\Delta>1$. We
can write $L^{(N)}=L_1^{(N)}+L_2^{(N)}+L_3^{(N)}$ depending on whether
$|\Psi(X_n)|>\Delta N$, $|\Psi(X_n)|\in(\Delta^{-1} N,\Delta N]$,
or $|\Psi(X_n)|\le(\Delta)^{-1} N$. Then
\begin{eqnarray*}
L_1^{(N)}\le\sum_{n=1}^{N}{\mathbb P}\bigl[|\Psi(X_n)|>\Delta N,
\bigl|Z_n^{(N)}\bigr|\le
N\bigr]=N{\mathbb P}\bigl[|\Psi(X_1)|>\Delta N, \bigl|Z_1^{(N)}\bigr|\le N\bigr].
\end{eqnarray*}
From formula (\ref{031023}) for $Z_n^{(N)}$, we conclude that the
event under the conditional probability can take place only when
$|P\chi_N(X_n)|$, or $|P\chi_N(X_{n-1})|>N(\Delta-1)/3$ for those
$N$, for which $c_N/N\le(\Delta-1)/3$. Using this observation, (\ref
{020411b}) and Chebyshev's inequality, one
can easily see that
\[
L_1^{(N)}\le2N[N(\Delta-1)/3]^{-\alpha'}\|P\chi_N\|_{L^{\alpha
'}(\pi
)}^{\alpha'}\to0
\]
as $N\to+\infty$. To deal with
$L_2^{(N)}$ consider a nonnegative $g\in C^{\infty}({\mathbb R})$
such that
$\|g\|_{\infty}\le1$, $[\Delta^{-1},\Delta]\subset\operatorname{supp}
g\subset[\Delta^{-1}_1,\Delta_1]$ for some $\Delta_1>\Delta$.
Repeating the foregoing argument for $g_2$, we conclude that
\[
\limsup_{N\to+\infty}L_2^{(N)}\le\|g\|_\infty\int_{\Delta
^{-1}_1}^{\Delta_1}\frac{d\lambda}{\lambda^2},
\]
which can be made as small as we wish by choosing $\Delta_1$
sufficiently close to $1$.
As for $L_3^{(N)}$, note that it equals
%
%
\begin{equation}\hspace*{29pt}
\label{041024}
L_3^{(N)}=\frac{1}{N}{\mathbb E} \Biggl|\sum_{n=1}^{N}{\mathbb E}
\bigl[M_n^{(N)}1\bigl[\bigl|M_n^{(N)}\bigr|\le N, |\Psi(X_n)|
\le(\Delta)^{-1}N\bigr]\big|\mathcal{G}_{n-1}\bigr] \Biggr|,
\end{equation}
where
%
%
\begin{eqnarray}
\label{012410}
M_n^{(N)}:\!&=&\chi_N(X_n)-P\chi_N(X_{n-1})
\nonumber\\[-8pt]\\[-8pt]
&=&\Psi_N(X_n)-c_N+P\chi
_N(X_{n})-P\chi_N(X_{n-1}).\nonumber
\end{eqnarray}
Thanks to the fact that $M_n^{(N)}$ are martingale differences, the
expression in (\ref{041024}) can be written as
$L_3^{(N)}=-(L_{31}^{(N)}+L_{32}^{(N)}+L_{33}^{(N)})$ where
$L_{3i}^{(N)}$ correspond to taking the conditional expectation over
the events $A_i$ for $i=1,2,3$
given by
\begin{eqnarray*}
A_1&:=&\bigl[\bigl|M_n^{(N)}\bigr|> N, |\Psi(X_n)|\le(\Delta)^{-1}N\bigr],\\
A_2&:=&\bigl[\bigl|M_n^{(N)}\bigr|> N, |\Psi(X_n)|>(\Delta)^{-1}N\bigr],\\
A_3&:=&\bigl[\bigl|M_n^{(N)}\bigr|\le N, |\Psi(X_n)|> (\Delta)^{-1}N\bigr].
\end{eqnarray*}
To estimate $L_{3i}^{(N)}$, $i=1,2$ we
note from (\ref{012410}) that $|M_n^{(N)}|> N$ only when $\Psi
_N(X_n)=\Psi(X_n)$ and $|\Psi(X_n)|\le N$, or $P\chi_N(X_{n-1})$,
$P\chi_N(X_n)$ are greater than $cN$ for some $c>0$. In the latter two
cases we can estimate similarly to $L_1^{(N)}$. In the first one,
however, we end up with the limit
\begin{eqnarray*}
&&
\limsup_{N\to+\infty}\frac{1}{N}{\mathbb E}\sum_{n=1}^{N}{\mathbb
E}
\bigl[ \bigl(|\Psi_N(X_n)|+|c_N|
+|P\chi_N(X_n)|+|P\chi_N(X_{n-1})| \bigr),\\
&&\hspace*{178.32pt}N\ge|\Psi
(X_n)|>\Delta^{-1}N |\mathcal{G}_{n-1} \bigr]\\
&&\qquad
\le\limsup_{N\to+\infty}N(1+|c_N|/N)\pi[N\ge|\Psi|> (\Delta)^{-1}N]
\\
&&\qquad\quad{}
+\limsup_{N\to+\infty}\int(I+P)|P\chi_N| 1[N\ge|\Psi|> (\Delta
)^{-1}N]\,d\pi.
\end{eqnarray*}
The second term on the utmost right-hand side vanishes thanks to (\ref
{020411b}).
The first one can be estimated as in the proof of Lemma \ref{lm:3},
and we obtain
that it is smaller than
$
C\int_{\Delta^{-1}}^1\lambda^{-2}\,d\lambda,
$
which can be made as small as we wish upon choosing $\Delta$
sufficiently close to $1$. We can estimate, therefore,
\begin{eqnarray*}
&&
\limsup_{N\to+\infty}L_{33}^{(N)}\\
&&\qquad
\le\limsup_{N\to+\infty}\frac{1}{N}{\mathbb E}\sum
_{n=1}^{N}{\mathbb E}
\bigl[ \bigl(|\Psi_N(X_n)|+|c_N| \bigr),N\ge|\Psi(X_n)|>\Delta
^{-1}N |\mathcal{G}_{n-1} \bigr]\\
&&\qquad\quad{}
+\limsup_{N\to+\infty}\frac{1}{N}{\mathbb E}\sum_{n=1}^{N}{\mathbb
E}
\bigl[ \bigl(|P\chi_N(X_n)|+|P\chi_N(X_{n-1})| \bigr),\\
&&\hspace*{158pt} |\Psi
(X_n)|>\Delta^{-1}N |\mathcal{G}_{n-1} \bigr]\\
&&\qquad
=\limsup_{N\to+\infty}
N\pi[N\ge|\Psi|> (\Delta)^{-1}N]\\
&&\qquad\quad{}+
\limsup_{N\to+\infty}\int(I+P)|P\chi_N| 1[ |\Psi|> (\Delta
)^{-1}N]\,d\pi
\\
&&\qquad
\le C\int_{\Delta^{-1}}^1\frac{d\lambda}{\lambda^2},
\end{eqnarray*}
which again can be made arbitrarily small.

\section[Proof of Theorem 2.8]{Proof of Theorem \protect\ref{thm-main-3}}
\label{cont-time}


Suppose that we are given a sequence of i.i.d. nonnegative random
variables $\{\rho_n, n\ge0\}$ independent of
$\{X_n, n\ge0\}$
and such that $A_{\alpha}:=\int_0^{+\infty}\rho^{\alpha}\varphi
(d\rho
)<+\infty$, where $\varphi(\cdot)$ is the distribuant of $\rho_0$
and $\alpha\in(0,2)$.
We consider a slightly more general situation than the one presented in
Theorem~\ref{stable-discrete-1} by allowing
%
%
\begin{equation}\label{eq:13}
S_N(t):=\sum_{n=0}^{[Nt]} \Psi(X_n)\rho_n.
\end{equation}
Observe that, if $\pi$ is the law of $X_n$, observable $\Psi$
satisfies the tail conditions
(\ref{eq:tails}), and $\rho_n$ is independent of $X_n$, then
\[
\lambda^\alpha\mathbb P \bigl( \Psi(X_0)\rho_0 > \lambda\bigr) =
\int_0^\infty
\lambda^\alpha\pi(\Psi> \lambda\rho^{-1}) \varphi(d\rho)
\mathop{\longrightarrow}_{\lambda\to+\infty}  c_*^+ A_\alpha.
\]
Define also
%
%
\begin{equation}
\label{031024}
C_N:=\int_{|\Psi|\le N} \Psi \,d\pi\,{\mathbb E}\rho_0.
\end{equation}

Consider then the Markov chain $\{(X_n, \rho_n), n\ge0\}$ on
$E\times
{\mathbb R}_+$. This Markov chain satifies all conditions used in the previous
sections, with stationary ergodic measure given by $\pi(dy)\otimes
\varphi(d\rho)$.
Then with the same arguments as used in Section~\ref{sec5} we get the
following.
\begin{theorem}
\label{thm:1}
\textup{(i)} Under the assumptions of the respective part \textup{(i)}, or \textup{(ii)} of Theorem
\ref{stable-discrete-1},
we have
$N^{-1/\alpha}S_N(\cdot)\stackrel{\mathrm{f.d.}}{\Rightarrow} Z(\cdot)$, as
$N\to+\infty$ where $\{Z(t), t\ge0\}$ is an $\alpha$-stable
process of type either type \textup{I}, or \textup{II} with the parameters of the
corresponding
L\'{e}vy measure [cf. (\ref{c-levy})]
given by
%
%
\begin{equation}
\label{c-levy-b}
C_*(\lambda):=\cases{
\alpha A_{\alpha}c_*^-,&\quad when $\lambda<0$,\cr
\alpha A_{\alpha}c_*^+,&\quad when $\lambda>0$.}
\end{equation}
Here $\stackrel{\mathrm{f.d.}}{\Rightarrow}$ denotes the convergence in the
sense of finite-dimensional distributions.

\textup{(ii)} In addition, under the assumptions of part \textup{(iii)} of Theorem \ref
{stable-discrete-1}
finite-dimensional distributions of
$N^{-1}S_N(t)-C_Nt$ converge weakly to those of $\{Z(t), t\ge0\}$, a
stable process of type \textup{III}. Here
$C_N$ is given by (\ref{031024}).
\end{theorem}
\begin{Remark*}
The results of the first part of the above theorem follow
under the conditions of Theorem \ref{thm-coupl},
by using the coupling argument of Section \ref{sec:coupling-approach}.
\end{Remark*}

Let us consider now the process $Y_N(t)$ defined by (\ref{eq:101}).
We only show that one-dimensional distributions of $Y_N(t)$ converge
weakly to the respective distribution of a suitable stable process $\{
Z(t), t\ge0\}$. The proof of convergence of
finite-dimensional distributions can be done in the same way.

Given $t>0$ define $n(t)$ as the positive integer, such that
\[
t_{n(t)}\le t<t_{n(t)+1},
\]
where $t_N$ is given by (\ref{t-N}).
Let
\begin{eqnarray*}
s(t)&:=&t/\bar t, \\
B_N(t)&:=&N^{-1/\alpha}\sum_{k=0}^{[Nt]}\Psi(X_k)\tau_k,\qquad t\ge0,
\end{eqnarray*}
where, as we recall, $\Psi(x):=V(x)t(x)$, $x\in E$ and $\{\tau_k,
k\ge0\}$
is a sequence of i.i.d. variables distributed according to an
exponential distribution with parameter~$1$.
Using the ergodic theorem one can easily conclude that
%
%
\begin{equation}
\label{010203}
s_N(t):=\frac{n(Nt)}{N}\to s(t)\qquad \mbox{as }N\to+\infty,
\end{equation}
a.s. uniformly on intervals of the form $[t_0,T]$ where $0<t_0<T$.
%
We have
\begin{eqnarray*}
&&
Y_N(t)=\frac{1}{N^{1/\alpha}}\sum_{k=0}^{n(Nt)-1}\Psi(X_k)\tau_k
+\frac{Nt-t_{n(Nt)}}{N^{1/\alpha}}V(X_k).
\end{eqnarray*}
Note that
\[
\frac{1}{N^{1/\alpha}}\sum_{k=0}^{n(Nt)}\Psi(X_k)\tau_k=B_N(s_N(t)).
\]
\begin{lemma}
For any $t>0$ and $\varepsilon>0$ fixed, we have
%
%
\begin{equation}
\label{approx}
\lim_{N\to+\infty}\mathbb P[|Y_N(t)-B_N(s_N(t))|>\varepsilon]
=0.
\end{equation}
\end{lemma}
\begin{pf}
Let $\sigma>0$ be arbitrary. We can write that
%
%
\begin{eqnarray}
\label{020203}
&&\mathbb P[|Y_N(t)-B_N(s_N(t))|>\varepsilon]\nonumber\\
&&\qquad\le\mathbb
P[|s_N(t)-s(t)|>\sigma]
\\
&&\qquad\quad{}+\mathbb P[|s_N(t)-s(t)|\le\sigma,
|Y_N(t)-B_N(s_N(t))|>\varepsilon
].\nonumber
\end{eqnarray}
The second term on the right-hand side can be estimated from above by
\begin{eqnarray*}
&&
\mathbb P\bigl[|s_N(t)-s(t)|\le\sigma, N^{-1/\alpha}\bigl|\Psi
\bigl(X_{n(Nt)}\bigr)\bigr|\tau
_{n(Nt)}>\varepsilon\bigr]\\
&&\qquad\le
\mathbb P\bigl[\sup\bigl\{|\Psi(X_{k})|\tau_{k}\dvtx k\in\bigl[\bigl(s(t)-\sigma\bigr)N,
\bigl(s(t)+\sigma\bigr)N\bigr]\bigr\}>N^{1/\alpha}\varepsilon\bigr].
\end{eqnarray*}
Using the stationarity of $\{|\Psi(X_{k})|\tau_{k}, k\ge0\}$ the
term on the right-hand side equals
\begin{eqnarray*}
&&
\mathbb P\bigl[\sup\{|\Psi(X_{k})|\tau_{k}\dvtx
k\in[0,2\sigma N]\}>N^{1/\alpha}\varepsilon\bigr]\\
&&\qquad\le2\sigma N\int_0^{+\infty}e^{-\tau}\pi[|\Psi(x)|\ge\tau
^{-1} N^{1/\alpha}\varepsilon]\,d\tau\le\frac{C\sigma}{\varepsilon
^{\alpha}}
\end{eqnarray*}
for some constant $C>0$, by virtue of (\ref{eq:tails}).
From (\ref{020203}) we obtain, therefore,
\begin{eqnarray*}
&&
\limsup_{N\to+\infty}\mathbb P[|Y_N(t)-B_N(s_N(t))|>\varepsilon]\le
\frac
{C\sigma}{\varepsilon^{\alpha}}
\end{eqnarray*}
for an arbitrary $\sigma>0$, which in turn implies (\ref{approx}).
\end{pf}

It suffices, therefore, to prove that the laws of $B_N(s_N(t))$
converge, as $N\to+\infty$, to the law of the respective stable process.
According to Skorochod's embedding theorem, one can find pairs of
random elements
$(\tilde B_N(\cdot),\tilde s_N(t))$, $N\ge1$, with values in
$D[0,+\infty)\times[0,+\infty)$, such that the law of each pair is
identical with that of $(B_N(\cdot), s_N(t))$, and $(\tilde B_N(\cdot
),\tilde s_N(t))$
converges a.s., as $N\to+\infty$, in the Skorochod topology to
$(Z(\cdot),s(t))$. Here $\{Z(t), t\ge0\}$ is the
stable process, as in Theorem \ref{thm:1}.
According to Proposition 3.5.3 page 119 of \cite{ethierkurtz},
the above means
that for each $T>0$ there exists a sequence of increasing
homeomorphisms $\lambda_N\dvtx[0,T]\to[0,T]$ such that
%
%
\begin{equation}
\label{010204}
\lim_{N\to+\infty}\gamma(\lambda_N)=0,
\end{equation}
where
\[
\gamma(\lambda_N):=\sup_{0<s<t<T} \biggl|\log\frac{\lambda
_N(t)-\lambda_N(s)}{t-s} \biggr|=0
\]
and
%
%
\begin{equation}
\label{030203}
{\sup_{t\in[0,T]}}|\tilde B_N\circ\lambda_N(t)-Z(t)|=0.
\end{equation}
As a consequence of (\ref{010204}) we have of course
that
%
%
\begin{equation}
\label{020204}
{\lim_{N\to+\infty}\sup_{t\in[0,T]}}|\lambda_N(t)-t|=0.
\end{equation}
Note that
the law of each $B_N(s_N(t))$ is identical with that of
$\tilde B_N(\tilde s_N(t))$. We also have
\begin{eqnarray*}
|\tilde B_N(\tilde s_N(t))-Z(s(t))|&\le&
|\tilde B_N(\tilde s_N(t))-Z\circ\lambda_N^{-1}(\tilde s_N(t))|
\\
&&{}+|Z\circ\lambda_N^{-1}(\tilde s_N(t))-Z(s(t))|.
\end{eqnarray*}
The right-hand side, however, vanishes a.s., as $N\to+\infty$, thanks
to (\ref{030203}), (\ref{020204}) and the fact that for each fixed
$s>0$ one has
$
{\mathbb P}[Z(s-)=Z(s)]=1
$
(see, e.g., Theorem 11.1, page 59 of \cite{sato}). The above allows us
to conclude that
$|\tilde B_N(\tilde s_N(t))-Z( s(t))|\to0$ a.s., as $N\to+\infty$,
thus the assertions of Theorem~\ref{thm-main-3} follow.

\section*{Acknowledgments}
The authors wish to express
their thanks
to the anonymous referee for thorough reading of the manuscript and
useful remarks that lead to the improvement of the presentation.
Milton Jara would like to
thank the hospitality of Universit\'{e} Paris-Dauphine and Maria
Curie-Sklodowska University (Lublin), where part of this work has been
accomplished.


%
\printaddresses


\begin{thebibliography}{27}

\bibitem{BPRNonlin}
\begin{barticle}[mr]
\bauthor{\bsnm{Bal},~\bfnm{Guillaume}\binits{G.}},
  \bauthor{\bsnm{Papanicolaou},~\bfnm{George}\binits{G.}} \AND
  \bauthor{\bsnm{Ryzhik},~\bfnm{Leonid}\binits{L.}}
(\byear{2002}).
\btitle{Radiative transport limit for the random {S}chr\"odinger equation}.
\bjournal{Nonlinearity}
\bvolume{15}
\bpages{513--529}.
\bid{mr={1888863}}
\end{barticle}
\endbibitem

\bibitem{basollaspohn}
\begin{barticle}[unstr]
\bauthor{\bsnm{Basile},~\bfnm{G.}\binits{G.}},
\bauthor{\bsnm{Olla},~\bfnm{S.}\binits{S.}} \AND
\bauthor{\bsnm{Spohn},~\bfnm{H.}\binits{H.}}
(\byear{2009}).
\btitle{Energy transport in
  stochastically perturbed lattice dynamics}.
\bjournal{Arch. Ration. Mech. Anal.}
DOI:
\href{http://dx.doi.org/10.1007/s00205-008-0205-6}{10.1007/s00205-008-0205-6}.
To appear.%
\end{barticle}%
\endbibitem%


\bibitem{browneagleson}
\begin{barticle}[mr]
\bauthor{\bsnm{Brown},~\bfnm{B.~M.}\binits{B.~M.}} \AND
  \bauthor{\bsnm{Eagleson},~\bfnm{G.~K.}\binits{G.~K.}}
(\byear{1971}).
\btitle{Martingale convergence to infinitely divisible laws with finite
  variances}.
\bjournal{Trans. Amer. Math. Soc.}
\bvolume{162}
\bpages{449--453}.
\bid{mr={0288806}}
\end{barticle}
\endbibitem

\bibitem{cc}
\begin{barticle}[mr]
\bauthor{\bsnm{Cs{\'a}ki},~\bfnm{Endre}\binits{E.}} \AND
  \bauthor{\bsnm{Cs{\"o}rg{\"o}},~\bfnm{Mikl{\'o}s}\binits{M.}}
(\byear{1995}).
\btitle{On additive functionals of {M}arkov chains}.
\bjournal{J. Theoret. Probab.}
\bvolume{8}
\bpages{905--919}.
\bid{mr={1353559}}
\end{barticle}
\endbibitem

\bibitem{demassi}
\begin{barticle}[mr]
\bauthor{\bsnm{De~Masi},~\bfnm{A.}\binits{A.}},
  \bauthor{\bsnm{Ferrari},~\bfnm{P.~A.}\binits{P.~A.}},
  \bauthor{\bsnm{Goldstein},~\bfnm{S.}\binits{S.}} \AND
  \bauthor{\bsnm{Wick},~\bfnm{W.~D.}\binits{W.~D.}}
(\byear{1989}).
\btitle{An invariance principle for reversible {M}arkov processes.
  {A}pplications to random motions in random environments}.
\bjournal{J. Stat. Phys.}
\bvolume{55}
\bpages{787--855}.
\bid{mr={1003538}}
\end{barticle}
\endbibitem

\bibitem{derrenniclin}
\begin{barticle}[mr]
\bauthor{\bsnm{Derriennic},~\bfnm{Yves}\binits{Y.}} \AND
  \bauthor{\bsnm{Lin},~\bfnm{Michael}\binits{M.}}
(\byear{2003}).
\btitle{The central limit theorem for {M}arkov chains started at a point}.
\bjournal{Probab. Theory Related Fields}
\bvolume{125}
\bpages{73--76}.
\bid{mr={1952457}}
\end{barticle}
\endbibitem

\bibitem{doeblin}
\begin{barticle}[mr]
\bauthor{\bsnm{Doeblin},~\bfnm{Wolfgang}\binits{W.}}
(\byear{1938}).
\btitle{Sur deux probl\`emes de {M}. {K}olmogoroff concernant les cha\^{\i}nes
  d\'enombrables}.
\bjournal{Bull. Soc. Math. France}
\bvolume{66}
\bpages{210--220}.
\bid{mr={1505091}}
\end{barticle}
\endbibitem

\bibitem{durrett}
\begin{bbook}[mr]
\bauthor{\bsnm{Durrett},~\bfnm{Richard}\binits{R.}}
(\byear{1996}).
\btitle{Probability: Theory and Examples}, \bedition{2nd} ed.
\bpublisher{Duxbury Press}, \baddress{Belmont, CA}.
\bid{mr={1609153}}
\end{bbook}
\endbibitem

\bibitem{durrettresnick}
\begin{barticle}[mr]
\bauthor{\bsnm{Durrett},~\bfnm{Richard}\binits{R.}} \AND
  \bauthor{\bsnm{Resnick},~\bfnm{Sidney~I.}\binits{S.~I.}}
(\byear{1978}).
\btitle{Functional limit theorems for dependent variables}.
\bjournal{Ann. Probab.}
\bvolume{6}
\bpages{829--846}.
\bid{mr={503954}}
\end{barticle}
\endbibitem

\bibitem{ErdosYau}
\begin{barticle}[mr]
\bauthor{\bsnm{Erd{\"o}s},~\bfnm{L{\'a}szl{\'o}}\binits{L.}} \AND
  \bauthor{\bsnm{Yau},~\bfnm{Horng-Tzer}\binits{H.-T.}}
(\byear{2000}).
\btitle{Linear {B}oltzmann equation as the weak coupling limit of a random
  {S}chr\"odinger equation}.
\bjournal{Comm. Pure Appl. Math.}
\bvolume{53}
\bpages{667--735}.
\bid{mr={1744001}}
\end{barticle}
\endbibitem

\bibitem{ethierkurtz}
\begin{bbook}[mr]
\bauthor{\bsnm{Ethier},~\bfnm{Stewart~N.}\binits{S.~N.}} \AND
  \bauthor{\bsnm{Kurtz},~\bfnm{Thomas~G.}\binits{T.~G.}}
(\byear{1986}).
\btitle{Markov Processes: Characterization and Convergence}.
\bpublisher{Wiley}, \baddress{New York}.
\bid{mr={838085}}
\end{bbook}
\endbibitem

\bibitem{fannjiang}
\begin{barticle}[mr]
\bauthor{\bsnm{Fannjiang},~\bfnm{Albert~C.}\binits{A.~C.}}
(\byear{2005}).
\btitle{White-noise and geometrical optics limits of {W}igner--{M}oyal equation
  for wave beams in turbulent media}.
\bjournal{Comm. Math. Phys.}
\bvolume{254}
\bpages{289--322}.
\bid{mr={2117627}}%
\end{barticle}%
\endbibitem%

\bibitem{foque}
\begin{bbook}[mr]
\bauthor{\bsnm{Fouque},~\bfnm{Jean-Pierre}\binits{J.-P.}},
  \bauthor{\bsnm{Garnier},~\bfnm{Josselin}\binits{J.}},
  \bauthor{\bsnm{Papanicolaou},~\bfnm{George}\binits{G.}} \AND
  \bauthor{\bsnm{S{\o}lna},~\bfnm{Knut}\binits{K.}}
(\byear{2007}).
\btitle{Wave Propagation and Time Reversal in Randomly Layered Media}.
\bseries{Stochastic Modelling and Applied Probability}
\bvolume{56}.
\bpublisher{Springer}, \baddress{New York}.
\bid{mr={2327824}}
\end{bbook}
\endbibitem

\bibitem{goldstein}
\begin{barticle}[mr]
\bauthor{\bsnm{Goldstein},~\bfnm{Sheldon}\binits{S.}}
(\byear{1995}).
\btitle{Antisymmetric functionals of reversible {M}arkov processes}.
\bjournal{Ann. Inst. H. Poincar\'e Probab. Statist.}
\bvolume{31}
\bpages{177--190}.
\bid{mr={1340036}}
\end{barticle}
\endbibitem

\bibitem{gordin}
\begin{barticle}[mr]
\bauthor{\bsnm{Gordin},~\bfnm{M.~I.}\binits{M.~I.}}
(\byear{1969}).
\btitle{The central limit theorem for stationary processes}.
\bjournal{Dokl. Akad. Nauk}
\bvolume{188}
\bpages{739--741}.
\bid{mr={0251785}}
\end{barticle}
\endbibitem

\bibitem{v1}
\begin{bmisc}[unstr]
\bauthor{\bsnm{Jara},~\bfnm{M.}\binits{M.}},
\bauthor{\bsnm{Komorowski},~\bfnm{T.}\binits{T.}} \AND
\bauthor{\bsnm{Olla},~\bfnm{S.}\binits{S.}}
(\byear{2009}).
\bhowpublished{Limit theorems for
additive functionals of a Markov Chain. Available at}
\href{http://arxiv.org/abs/0809.0177v1}{arXiv:0809.0177v1}.
\end{bmisc}
\endbibitem

\bibitem{kipnislandim}
\begin{bbook}[mr]
\bauthor{\bsnm{Kipnis},~\bfnm{Claude}\binits{C.}} \AND
  \bauthor{\bsnm{Landim},~\bfnm{Claudio}\binits{C.}}
(\byear{1999}).
\btitle{Scaling Limits of Interacting Particle Systems}.
\bseries{Grundlehren der Mathematischen Wissenschaften [Fundamental Principles
  of Mathematical Sciences]}
\bvolume{320}.
\bpublisher{Springer}, \baddress{Berlin}.
\bid{mr={1707314}}
\end{bbook}
\endbibitem

\bibitem{kipnisvaradhan}
\begin{barticle}[mr]
\bauthor{\bsnm{Kipnis},~\bfnm{C.}\binits{C.}} \AND
  \bauthor{\bsnm{Varadhan},~\bfnm{S.~R.~S.}\binits{S.~R.~S.}}
(\byear{1986}).
\btitle{Central limit theorem for additive functionals of reversible {M}arkov
  processes and applications to simple exclusions}.
\bjournal{Comm. Math. Phys.}
\bvolume{104}
\bpages{1--19}.
\bid{mr={834478}}
\end{barticle}
\endbibitem

\bibitem{lax}
\begin{bbook}[mr]
\bauthor{\bsnm{Lax},~\bfnm{Peter~D.}\binits{P.~D.}}
(\byear{2002}).
\btitle{Functional Analysis}.
\bpublisher{Wiley}, \baddress{New York}.
\bid{mr={1892228}}
\end{bbook}
\endbibitem

\bibitem{sll}
\begin{barticle}[mr]
\bauthor{\bsnm{Lepri},~\bfnm{Stefano}\binits{S.}},
  \bauthor{\bsnm{Livi},~\bfnm{Roberto}\binits{R.}} \AND
  \bauthor{\bsnm{Politi},~\bfnm{Antonio}\binits{A.}}
(\byear{2003}).
\btitle{Thermal conduction in classical low-dimensional lattices}.
\bjournal{Phys. Rep.}
\bvolume{377}
\bpages{1--80}.
\bid{mr={1978992}}
\end{barticle}
\endbibitem

\bibitem{LukSpohn}
\begin{barticle}[mr]
\bauthor{\bsnm{Lukkarinen},~\bfnm{Jani}\binits{J.}} \AND
  \bauthor{\bsnm{Spohn},~\bfnm{Herbert}\binits{H.}}
(\byear{2007}).
\btitle{Kinetic limit for wave propagation in a random medium}.
\bjournal{Arch. Ration. Mech. Anal.}
\bvolume{183}
\bpages{93--162}.
\bid{mr={2259341}}
\end{barticle}
\endbibitem

\bibitem{maxwellwoodrofe}
\begin{barticle}[mr]
\bauthor{\bsnm{Maxwell},~\bfnm{Michael}\binits{M.}} \AND
  \bauthor{\bsnm{Woodroofe},~\bfnm{Michael}\binits{M.}}
(\byear{2000}).
\btitle{Central limit theorems for additive functionals of {M}arkov chains}.
\bjournal{Ann. Probab.}
\bvolume{28}
\bpages{713--724}.
\bid{mr={1782272}}
\end{barticle}
\endbibitem

\bibitem{mmm}
\begin{bmisc}[unstr]
\bauthor{\bsnm{Mellet},~\bfnm{A.}\binits{A.}},
\bauthor{\bsnm{Mischler},~\bfnm{S.}\binits{S.}} \AND
\bauthor{\bsnm{Mouhot},~\bfnm{C.}\binits{C.}}
(\byear{2009}).
\bhowpublished{Fractional diffusion
limit for collisional kinetic equations. Available at}
\href{http://arxiv.org/abs/0809.2455}{arXiv:0809.2455}.
\end{bmisc}
\endbibitem

\bibitem{nagaev}
\begin{barticle}[mr]
\bauthor{\bsnm{Nagaev},~\bfnm{S.~V.}\binits{S.~V.}}
(\byear{1957}).
\btitle{Some limit theorems for stationary {M}arkov chains}.
\bjournal{Theory Probab. Appl.}
\bvolume{2}
\bpages{378--406}.
\bid{mr={0094846}}
\end{barticle}
\endbibitem

\bibitem{olla}
\begin{bincollection}[mr]
\bauthor{\bsnm{Olla},~\bfnm{Stefano}\binits{S.}}
(\byear{2001}).
\btitle{Central limit theorems for tagged particles and for diffusions in
  random environment}.
In \bbooktitle{Milieux Al\'eatoires}.
\bseries{Panoramas et Synth\`eses}
\bvolume{12}
\bpages{23--25}.
\bpublisher{Soc. Math. France}, \baddress{Paris}.
\bid{mr={2226846}}
\end{bincollection}
\endbibitem

\bibitem{sato}
\begin{bbook}[mr]
\bauthor{\bsnm{Sato},~\bfnm{Ken}\binits{K.}}
(\byear{1999}).
\btitle{L\'evy Processes and Infinitely Divisible Distributions}.
\bseries{Cambridge Studies in Advanced Mathematics}
\bvolume{68}.
\bpublisher{Cambridge Univ. Press}, \baddress{Cambridge}.
\bid{mr={1739520}}
\end{bbook}
\endbibitem

\bibitem{spohn}
\begin{barticle}[mr]
\bauthor{\bsnm{Spohn},~\bfnm{Herbert}\binits{H.}}
(\byear{1977}).
\btitle{Derivation of the transport equation for electrons moving through
  random impurities}.
\bjournal{J. Stat. Phys.}
\bvolume{17}
\bpages{385--412}.
\bid{mr={0471824}}
\end{barticle}
\endbibitem

\end{thebibliography}
\end{document}